\documentclass[english]{amsart}
\usepackage[T1]{fontenc}
\usepackage[latin1]{inputenc}
\usepackage{a4}
\usepackage{babel}
\usepackage{graphics}

\makeatletter

\providecommand{\LyX}{L\kern-.1667em\lower.25em\hbox{Y}\kern-.125emX\@}
\newcommand{\noun}[1]{\textsc{#1}}

 \theoremstyle{plain}    
 \newtheorem{thm}{Theorem}[section]
 \numberwithin{equation}{section} %% Comment out for sequentially-numbered
 \numberwithin{figure}{section} %% Comment out for sequentially-numbered
 \theoremstyle{plain}    
 \newtheorem*{thm*}{Theorem} 
 \theoremstyle{plain}    
 \newtheorem{cor}[thm]{Corollary} %%Delete [thm] to re-start numbering
 \theoremstyle{plain}    
 \newtheorem{lem}[thm]{Lemma} %%Delete [thm] to re-start numbering
 \theoremstyle{plain}    
 \newtheorem{prop}[thm]{Proposition} %%Delete [thm] to re-start numbering
 \theoremstyle{remark}
 \newtheorem*{rem*}{Remark}
 \newenvironment{lyxlist}[1]
   {\begin{list}{}
     {\settowidth{\labelwidth}{#1}
      \setlength{\leftmargin}{\labelwidth}
      \addtolength{\leftmargin}{\labelsep}
      }}
   {\end{list}}

\usepackage[T1]{fontenc}
\usepackage[latin1]{inputenc}
\usepackage{babel}

\makeatletter

\usepackage[mathscr]{eucal}
\usepackage{amssymb}
\makeatother

\makeatother
\AtBeginDocument{
  
}

\begin{document}
\newenvironment{proofing}{\medskip\textsc{Proof.}}{\hfill$\Box$\medskip} \newcommand{\NN}{\mathbb{N}} \newcommand{\RR}{\mathbb{R}} \newcommand{\ZZ}{\mathbb{Z}} \newcommand{\CC}{\mathbb{C}} \newcommand{\QQ}{\mathbb{Q}} \newcommand{\TT}{\mathbb{T}} \newcommand{\PP}{\mathbb{P}} \newcommand{\EE}{\mathbb{E}} \newcommand{\HH}{\mathbb{H}} \newcommand{\sA}{{\mathscr{A}}} \newcommand{\sB}{{\mathscr{B}}} \newcommand{\sE}{{\mathscr{E}}} \newcommand{\sF}{{\mathscr{F}}} \newcommand{\sG}{{\mathscr{G}}} \newcommand{\sH}{{\mathscr{H}}} \newcommand{\sJ}{{\mathscr{J}}} \newcommand{\sK}{{\mathscr{K}}} \newcommand{\sL}{{\mathscr{L}}} \newcommand{\sT}{{\mathscr{T}}} \newcommand{\sY}{{\mathscr{Y}}} \newcommand{\sZ}{{\mathscr{Z}}} \newcommand{\fN}{\mathfrak{N}} \newcommand{\fA}{\mathfrak{A}} \newcommand{\fS}{\mathfrak{S}} \newcommand{\fZ}{\mathfrak{Z}} \newcommand{\fp}{\mathfrak{p}} \newcommand{\fu}{\mathfrak{u}} \newcommand{\fs}{\mathfrak{s}} \newcommand{\esp}[2][]{\mathbb{E}_{#1}\!\left[#2\right]} \newcommand{\pr}[2][]{\mathbb{P}_{#1}\!\left[#2\right]} \newcommand{\id}{\mathrm{Id}} \newcommand{\ds}{\displaystyle} \newcommand{\limSup}[1][n]{\overline{\lim_{#1\to\infty}}} \newcommand{\iso}{\textrm{Iso}(\TT)}
\newcommand{\Rp}{\RR ^{*}_{+}}
 
\newcommand{\affR}{\mathrm{Aff}(\RR ^{d})}

\newcommand{\affT}{\mathrm{Aff}(\TT )}

\newcommand{\RxR}{\Rp \times \RR }

\newcommand{\affO}[1]{\mathrm{Aff}(#1 )}
 
\newcommand{\affq}{\mathrm{Aff}(\QQ _{p})}

\newcommand{\hor}{\mathrm{Hor}(\TT )}

\newcommand{\bT}{\partial ^{*}\TT }
 
\newcommand{\bG}{\partial ^{*}\Gamma }

\newcommand{\horg}{\mathrm{Hor}(\Gamma )}
 
\newcommand{\horO}[1]{\mathrm{Hor}(#1 )}

\newcommand{\toinf}{\rightarrow \infty }
 
\newcommand{\inv}{^{-1}}
 
\newcommand{\Ind}[1]{1_{#1 }}
 
\newcommand{\nor}[1]{\left| #1 \right| }
 
\newcommand{\sumzi}[1]{\sum ^{\infty }_{#1 =0}}
 
\newcommand{\dnor}[1]{\left\Vert #1 \right\Vert }

\newcommand{\tm}{\frac{3}{2}}

\newcommand{\stard}{\stackrel{\cdot }{*}}

\newcommand{\supp}{\mathrm{supp}}

\newcommand{\cf}{\wedge }

\newcommand{\e}{\mathrm{e}}

\newcommand{\ver}{\rightarrow }

\newcommand{\Qp}{\QQ _{p}}

{\raggedleft \today{}\par}
\medskip{}

\title{Renewal theory on the oriented tree }

\maketitle
\medskip{}
{\centering \noun{Sara Brofferio}%
\footnote{Institut für Mathematik C, Technische Universtät Graz, Austria

e-mail: brofferio@finanz.math.tu-graz.ac.at

This work has been partially supported by the Austrian Science Fund
(FWF), Project No. P15577-N05
}\par}
\bigskip{}

\begin{abstract}
The affine group of a tree is the group of the isometries of a homogeneous
tree that fix an end of its boundary. Consider a probability measure
\( \mu  \) on this group and the associated random walk. The main
goal of this paper is to determine the accumulation points of the
potential kernel \[
g*U=g*\sum _{n=0}^{\infty }\mu ^{(n)}\]
when \( g \) tends to infinity. In particular we show that under
suitable regularity hypotheses this kernel can be continuously extended
to the tree boundary and we determine the limit measures.

\emph{Key words:} Random walk, renewal theory, affine group, tree,
\( p \)-adic rationals 
\end{abstract}
\tableofcontents{}

\section*{Introduction}

Consider a transient random walk with law \( \mu  \) on a locally
compact group. Its potential measure \( U=\sum _{n=0}^{\infty }\mu ^{(n)} \)
is a Radon measure and its (right) potential kernel \( g*U \) is,
when \( g \) varies on the group, a family of measures that is vaguely
relatively compact. Renewal theory consists in studying the limit
behavior of this family when \( g \) goes to infinity, determining
the limit measures and the geometrical directions along which it converges.
On Abelian groups, this problem has been completely solved (cf. \cite{PS69}):
there are not more then two accumulation points (the null measure
and the Haar measure) and there is a non-zero limit if and only if
the group is a compact extension of \( \ZZ  \) or \( \RR  \). The
work of L.Elie on the affine group of the real line and on almost
connected Lie groups (\cite{El82}) has shown that for non-unimodular
groups we may have a quite different behavior. Namely there exists
an infinite number of limit measures, and the Haar measure cannot
be among them.

In this paper we leave the Euclidean setting to deal with this kind
of question on the group of affine transformations of the homogeneous
tree, \( \affT  \), that is the group of tree isometries that fix
an end of the boundary. D.Cartwright, V.Kaimanovich and W.Woess have
given in \cite{CKW} a first detailed study of random walks on this
group and we refer to this article for a comprehensive introduction.

The affine group of the tree contains the affine group of the \( p \)-adic
numbers, \( \affO{\Qp } \) (or more generally of a local field),
that is the group of matrices of the form \( \left[ \begin{array}{cc}
a & b\\
0 & 1
\end{array}\right]  \) where the coefficients \( a\neq 0 \) and \( b \) are \( p \)-adic
valued. The tree is in fact the Bruhat-Tits building of the invertible
\( 2\times 2 \) matrices on \( \Qp  \) and its affine group acts
on the tree analogously as the real affine group, \( \affO{\RR } \),
acts on the hyperbolic plane \( \HH ^{2} \), that is by isometries
and fixing a boundary point. On the other hand, the structural analogies
apart, the real affine group and the affine group of the tree present
remarkable differences. While \( \affO{\RR } \) can be identified
with the group of all isometries that fix a boundary point, the affine
group of the tree is significantly bigger and more complex then \( \affO{\Qp } \)
and it contains other interesting subgroups such the lamplighter group
or automatic groups. This complexity is mainly due to the fact that
the graph structure of the tree is much less rigid then the hyperbolic
plane, in the sense that the local behavior of an isometry does not
determine how it acts globally.

The main goal of this paper is to show that the potential kernel of
a random walk supported by non-exceptional subgroup of \( \affT  \)
can be continuously extended to the boundary of the tree and to give
a description of the limit measures by mean of the invariant measure
on the boundary and the counting measure on the integers (Theorems
\ref{teo-lim.U.in.bT} and \ref{teo-lim.pot.omega}). These last conclusions
are particularly interesting in view of farther studies, namely for
the characterization of the Martin boundary points and thus in the
representation of the invariant measures. 

Our results are obtained partially by adapting Elie's methods that
involve the characterization of the periods of the limit measures,
partially by using a weighty renewal equality (Corollary \ref{cor-eq.ren.AffT})
whose analogue over \( \RR  \) is due to M.Babillot, Ph.Bougerol
and L.Elie \cite{BBE}. In a general setting we require, besides weak
moment conditions, that the the step law of the random walk is spread
out. However, for random walks supported by groups that act on the
tree on a sufficiently homogeneous way ( such as \( \affO{\Qp } \)
) we are able to to avoid this last continuity hypothesis for the
limit toward all boundary points except for the one that is fixed
by the group. 

The paper is structured as follow:

In Section 1, we introduce the structures we are working on (the oriented
tree, the affine group and its non-exceptional subgroups) and the
probabilistic objects we are going to study (random walks and potential
kernel).

In Section 2, we give some preliminary results concerning the convergence
and the action of the random walks on the tree boundary and obtain
a measure equality for the potential kernel on the group.

In Section 3, we prove our main results. We start by determining some
invariance properties of the limit measures and then we characterize
the limits of the potential kernel.

\section{Random walks on the affine group of a tree}

\subsection{Oriented tree}

We consider the \emph{homogeneous tree} \( \TT  \) of degree \( q+1 \),
i.e. the connected non-oriented graph without cycles whose vertices
have exactly \( q+1 \) neighbors, equipped with the usual graph distance
\[
d(x,y)=\textrm{number of edges between }x\textrm{ and }y.\]
 The set of infinite geodesic rays that start from some vertex and
go to infinity, quotiented by the equivalence relation that identifies
two geodesics when they coincide but for a finite number of vertices
give the geometrical \emph{boundary} of the tree, \( \partial \TT  \).
The union \( \TT \cup \partial \TT  \), equipped with the topology
of the infinite cones starting from a vertex, is then a compact set
where \( \TT  \) is a dense open sub-set. 

A partial order of the tree is given, fixing an end \( \omega  \)
in \( \partial \TT  \) and setting for all \( x\neq y \) in \( \TT \cup \partial \TT  \)
\[
x\cf y=\textrm{ first common vertex of }\overline{x\omega }\textrm{ and }\overline{y\omega },\]
where \( \overline{x\omega } \) is the geodesic starting at \( x \)
and in the class of \( \omega  \), and \( x\wedge x=x \). We write
\[
x\succeq y\quad \Leftrightarrow \quad x=x\cf y.\]
 One can imagine the oriented tree as an infinite genealogical tree,
where \( \omega  \) represents the \emph{mythical ancestor}, every
vertex has \( q \) sons and a father and \( x\succeq y \) if and
only if \( y \) is a descendent of \( x \). 

\vspace{0.3cm}
{\centering \resizebox*{0.8\textwidth}{!}{\includegraphics{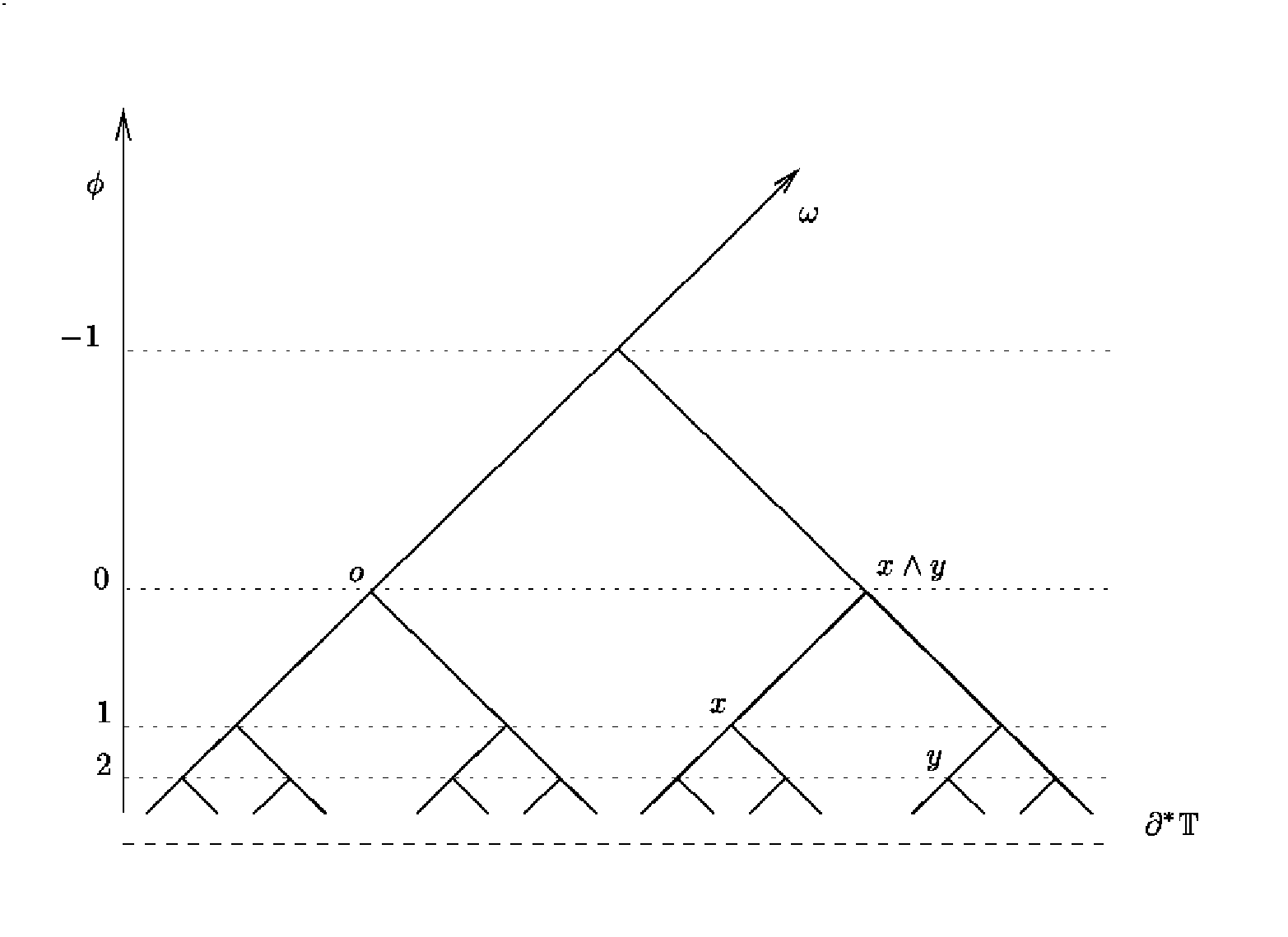}} \par}
\vspace{0.3cm}

Let us fix a reference vertex \( o \) in \( \TT  \) called \emph{origin}.
The height function \( \phi  \) from \( \TT  \) on \( \ZZ  \) is
\[
\phi (x):=d(x,x\cf o)-d(o,x\cf o),\]
also known as the \emph{Busemann} function, represents the generation
number of \( x \). 

Let consider the bottom boundary of the tree \[
\bT =\partial \TT -\{\omega \}.\]
The function \( \phi  \) induces a ultra-metric distance on \( \TT \cup \bT  \)
defined by \[
\Theta (\alpha ,\beta ):=\left\{ \begin{array}{cc}
q^{-\phi (\alpha \wedge \beta )}\textrm{ } & \textrm{if }\alpha \not =\beta \in \TT \cup \bT \\
0 & \textrm{if }\alpha =\beta 
\end{array}\right. .\]

\subsection{The affine group of the tree}

The group of isometries of the tree \( (\TT ,d) \) has a natural
continuous action on the boundary, obtained by the action on the geodesics. 

The \emph{affine group of the tree} is the subgroup of the isometries
that fix the end \( \omega  \) \[
\affT :=\left\{ g\in \iso :g\omega =\omega \right\} ,\]
 that is the subgroup that preserves the order induced by \( \omega  \),\[
g(x\cf y)=gx\cf gy\qquad \textrm{ for all }g\in \affT .\]

The group \( \affT  \) is equipped with the topology of pointwise
convergence on the tree, where a base of open neighborhoods of an
affinity \( \gamma  \) is given by the sets \[
V(\mathbf{x}\ver \mathbf{y}):=\left\{ g\in \affT :\, g\mathbf{x}=\mathbf{y}\right\} \]
for every finite set of vertices \( \mathbf{x}=(x_{1},\ldots ,x_{n}) \)
and with \( \mathbf{y}=\gamma \mathbf{x}=(\gamma x_{1},\ldots ,\gamma x_{n}) \).
These sets are simultaneously open and compact, therefore \( \affT  \)
is a locally compact totally disconnected group. 

The semi-norm \[
\nor{g}=d(go,o)\]
on \( \affT  \) is symmetric, \( \nor{g}=\nor{g^{-1}} \), and verifies
to \( \nor{g_{1}g_{2}}\leq \nor{g_{1}}+\nor{g_{2}} \). The set of
the affinities of zero norm is \( V(o\ver o) \), a compact subgroup.

\subsubsection{Drift of an affinity and the horocyclic group}

As the affinities respect the order  and the distance on the tree,
for every couple of vertices \( x \) and \( y \) one has\[
\phi (gx)-\phi (gy)=\phi (x)-\phi (y).\]
The homomorphism:\begin{eqnarray*}
\phi :\affT  & \rightarrow  & \ZZ \\
g & \mapsto  & \phi (gx)-\phi (x)=\phi (go),
\end{eqnarray*}
does not depend on the choice of the point \( x \) and contains the
information on vertical action of an affinity on the tree. It also
indicates whether the action of \( g \) on the bottom boundary \( \bT  \)
dilates or contracts, in fact for every couple of ends \( \alpha  \)
and \( \beta  \) in \( \bT  \)\begin{equation}
\label{eq-Theta(g)-phi(g)}
\Theta (g\alpha ,g\beta )=q^{-\phi (g\alpha \cf g\beta )}=q^{-\phi (g)}\Theta (\alpha ,\beta ).
\end{equation}

The \emph{horocyclic group} of the tree is the subgroup of the affine
group that fixes the heights \[
\hor :=\ker \phi =\left\{ g\in \affT :\phi (gx)=\phi (x)\quad \forall x\in \TT \right\} .\]
 It follows from (\ref{eq-Theta(g)-phi(g)}) that \( \hor  \) is
the group of all isometries of \( (\bT ,\Theta ) \). 

Instead of working on the whole affine group, we will be often interested
in some closed subgroup \( \Gamma  \) of \( \affT  \). In this case
we consider \( \phi  \) as an homomorphism from \( \Gamma  \) to
\( \ZZ  \) and set \[
\horg :=\ker \phi =\hor \cap \Gamma .\]

\subsubsection{Algebraic structure of \protect\( \affT \protect \)}

For the sake of simplicity we always suppose that the homomorphism
\( \phi  \) from \( \Gamma  \) on \( \ZZ  \) is surjective. Then
for all \( s\in \Gamma  \) such that \( \phi (s)=1 \), every \( g\in \Gamma  \)
has a unique decomposition as a product of an element of the horocyclic
group and a power of \( s \) \begin{equation}
\label{eq-g-bs}
g=b(g)s^{\phi (g)}\quad \textrm{ where }\quad b(g):=gs^{-\phi (g)}\in \horg .
\end{equation}
 Thus, if we identify \( \ZZ  \) with the subgroup generated by \( s \),
the group \( \Gamma  \) is the semi-direct product \begin{eqnarray*}
\horg \rtimes _{s}\ZZ  & \cong  & \Gamma \\
(b,h) & \mapsto  & bs^{h}.
\end{eqnarray*}

Note that the decomposition of \( \affT  \) as semi-direct product
of \( \ZZ  \) and \( \hor  \) depends on the choice of the element
\( s \), and we call it \emph{reference homothety}. We denote by
\( \alpha =\alpha _{s} \) the unique end of \( \bT  \) fixed by
\( s \) (for its existence see for instance \cite{Ti70}). The homothety
\( s \) acts by translation on the geodesic \( \overline{\alpha \omega } \),
that may then be considered as a {}``main branch'' of the tree.
To choose a reference homothety is equivalent to select a center \( \alpha  \)
of the bottom boundary and a canonical identification between the
sub-trees that branch from \( \overline{\alpha \omega } \).

\subsubsection{Rotations}

In some sense the horocyclic group, that is the group of all isometries
of the bottom boundary, plays the role of the group of translations
in the real case; but in our case the action on \( \bT  \) is not
simple. In fact the stabilizer of an end \( \alpha \in \bT  \) in
the horocyclic group , that is the group of \emph{rotations} of center
\( \alpha  \), is not trivial and is the compact subgroup \[
K_{\alpha }=K_{\alpha }(\Gamma )=\left\{ r\in \horg :\: r\alpha =\alpha \right\} .\]
 It is worthwhile observing that, contrary to what happens on Lie
groups, the identification induced by the reference homothety \( s \)
is completely arbitrary, because the structure of the tree is much
less rigid then in the analogous continuous spaces. One of the first
consequence is that a rotation does not commute with an homothety
of same center; in fact whenever a rotation \( r \) acts on two sub-trees
in different ways (according to the identification induced by \( s \)),
one has \( sr\neq rs \). Moreover as there is no rigidity, we do
not have a finite set of points whose images uniquely determine a
rotation or an affinity, but given any compact set \( C \) in \( \TT  \)
(or in \( \bT  \)) one may found two affinities that act in the same
way on \( C \) and are different on its complement.

\subsubsection{Compactification and group boundary}

The action on the tree enables us to give a natural compactification
of any closed subgroup of \( \affT  \). In fact it easy to see that
whenever for a sequence \( \left\{ g_{n}\right\} _{n} \) in \( \affT  \)
there exists a vertex \( x\in \TT  \) such that \( \left\{ g_{n}x\right\} _{n} \)
converges to an end \( \beta  \) in \( \partial \TT  \), then for
all \( y\in \TT  \) also \( \left\{ g_{n}y\right\} _{n} \) converges
to \( \beta  \). We then say that the sequence \( \left\{ g_{n}\right\} _{n} \)
\emph{converges to} \( \beta  \) and we compactify \( \affT  \)
in \( \affT \cup \partial \TT  \) setting \[
g_{n}\rightarrow \beta \in \partial \TT \Leftrightarrow \exists \, (\textrm{or }\forall )\, x\in \TT \, :\, \, g_{n}x\rightarrow \beta .\]

The boundary of a subgroup \( \Gamma  \) of \( \affT  \) is then
the set of the accumulation points of \( \Gamma  \) in \( \partial \TT  \)
and is denoted by \( \partial \Gamma  \).

\subsection{Non-exceptional subgroup}

We focus our study on random walks supported on subgroups of \( \affT  \)
that are non-degenerated. More precisely we deal with closed subgroups,
\( \Gamma  \), that are \emph{non-exceptional}, i.e. verify to one
of the following equivalent conditions (cf. \cite{CKW}): \emph{}

\begin{itemize}
\item \( \Gamma  \) is not contained in \( \hor  \) and it does not fix
any end in \( \bT  \)
\item \( \Gamma  \) is non-unimodular
\item \( \partial \Gamma  \) is infinite 
\end{itemize}
In the usual parallelism with the real affine group, this is equivalent
to ask that the group \( \Gamma  \) is neither a group of translations
nor of roto-homotheties. Another important property of non-exceptional
subgroups is that all their orbits are dense in the bottom boundary
of the group \[
\bG =\partial \Gamma -\left\{ \omega \right\} .\]
When the group is also closed then its action on \( \bG  \) is transitive,
i.e. for all \( \beta \in \bG  \)\[
\Gamma \beta =\bG .\]
 By (\ref{eq-g-bs}), also \( \horg  \) acts transitively on \( \bG  \).

Let \( s\in \Gamma  \) be a reference homothety of center \( \alpha  \).
As \( K_{\alpha }=K_{\alpha }(\Gamma ) \) is the stabilizer of the
end \( \alpha \in \bG  \) in \( \horg  \), the subgroup \( K_{\alpha }\rtimes _{s}\ZZ  \)
is the stabilizer of \( \alpha  \) in \( \Gamma  \); thus we have
the following identifications by homeomorphisms: \[
\Gamma /\left( K_{\alpha }\rtimes _{s}\ZZ \right) =\horg /K_{\alpha }=\bG .\]

\subsection{Random walks on  \protect\( \affT \protect \)}

Let \( \mu  \) be a probability measure on \( \affT  \) and \( \left\{ X_{n}\right\} _{n\in \NN } \)
a sequence of random variables defined on the probability space \( (\Omega ,\PP ) \)
and with values in \( \affT  \), independent and identically distributed
with law \( \mu  \). The \emph{left and right random walks} are the
Markov chains on \( \affT  \) defined by iterated products of the
\( X_{n} \) on the left and on the right, respectively, \[
L_{n}=X_{n}X_{n-1}\cdots X_{1}\qquad \textrm{and}\qquad R_{n}=X_{1}\cdots X_{n-1}X_{n}\]
and \( L_{0}=R_{0} \) are equal to the identity \( e \). 

Although from a trajectories point of view these two processes are
different, for every fixed time \( n \) have they the same law, the
\( n \)-th convolution power of \( \mu  \)\[
L_{n}\stackrel{\textrm{law}}{=}R_{n}\sim \mu ^{(n)}.\]

\subsubsection{Hypotheses}

All our results concern random walks whose action on the tree is sufficiently
complete, and namely we shall always assume the following \emph{non-degeneracy
hypotheses:} the closed subgroup generated by the support of \( \mu  \)\[
\Gamma :=\overline{<\mathrm{supp}\mu >}\]
 is non-exceptional and, just for sake of simplicity, that \( \phi (\Gamma )=\ZZ  \).

We also need some moment hypothesis that may vary according to the
type of results that we want to obtain. We always suppose that the
projection of the random walk on \( \ZZ  \) is integrable\[
\esp {\left| \phi (X_{1})\right| }<+\infty .\]
Most of the times we also require a \emph{moment of first order} for
the random walk on the group \[
\esp {\left| X_{1}\right| }<+\infty .\]
When projection of the random walk on \( \ZZ  \) is recurrent, we
shall need a \emph{moment of order \( 2+\epsilon  \),} namely \[
\esp {\phi (X_{1})^{2}+\nor{b(X_{1})}^{2+\varepsilon }}<+\infty \]
for some \( \varepsilon >0 \). 

For generic random walks on \( \affT  \) we also require a continuity
condition, namely that the measure \( \mu  \) is \emph{spread out}
(i.e. there exists a convolution power \( \mu ^{(n)} \) that is non-singular
with respect to the Haar measure of \( \Gamma  \)).

\subsubsection{Drift of the random walk}

A crucial role in the study of the random walks \( R_{n} \) and \( L_{n} \)
is played by their projection on \( \ZZ  \), that is by the random
walk \[
S_{n}=\phi (X_{1})+\cdots +\phi (X_{n})=\phi (L_{n})=\phi (R_{n}).\]
 Its mean is called \emph{drift} \emph{of} \( \mu  \) \[
\mu (\phi )=\esp {\phi (X_{1})}\]
 and it is the parameter that enables to classify the different types
of behavior.

\subsection{Remarkable examples}

Since the tree has a very lax structure, the affine group of the tree
is very complex and there are many different way to construct random
walks supported by non-exceptional subgroups. We present here some
remarkable examples.

\subsubsection{Random walks on the tree}

The simplest example is given by Markov chains on the tree that is
invariant under the transitive action of a subgroup \( \Gamma  \),
as for instance the nearest neighbor random walk on the tree where
from a vertex one goes to the father with probability \( \alpha  \)
and to every son with probability \( (1-\alpha )/q \). This type
of process can be obtained from a random walk on \( \Gamma  \), whose
law is invariant by the right action of the stabilizer \( V(x) \)
of a vertex \( x\in \TT  \). In that case the process \( Z_{n}=R_{n}V(x) \)
is a Markov chain on \( \Gamma /V(x)=\TT  \) starting in \( x \)
and homogeneous under the action of \( \Gamma  \). 

One can show that, since the stabilizer \( V(x) \) is an open and
compact subgroup, every measure that is right invariant by the action
of \( V(x) \) has a continuous density with respect to Haar measure.
Thus every Markov chain on the tree that is invariant under the transitive
action of an non-exceptional subgroup may be considered a random walk
on the group whose law is spread out, and therefore all our results
translated in this setting.

\subsubsection{\protect\( p\protect \)-adic affine group }

One of the most interesting non-exceptional subgroups of the affine
group of the tree is the affine group of rational \( p \)-adic. We
will often refer to it because, apart for its intrinsic interest,
it is the more natural generalization of the real affine group and
thus it allows to stress the similarities but also the main differences
between the real and tree settings.

Let \( p \) be a prime number, consider the integer evaluation \( v_{p} \)
on the rational numbers that measures how much a number is divisible
by \( p \), i.e. for every \( u\in \QQ ^{*} \), we set \[
v_{p}(u)=\max \{k\in \ZZ :\, p^{-k}ru\in \ZZ \}\]
and \( \nu _{p}(0)=0 \). The field of rational \( p \)-adic numbers
is then the completion of \( \QQ  \) equipped with the ultra-metric
norm\[
\nor{u}_{p}=p^{-v_{p}(u)}\quad \, \, \textrm{for all }u\in \QQ .\]

There exists a strict relationship between the \( p \)-adic rational
and the oriented tree of degree \( p+1 \) (cf. Serre \cite{Ser80}),
since it is possible to consider the tree as the set of the discs
\( \QQ _{p} \). First observe that, because the evaluation \( v_{p} \)
is integer valued and \( \nor{\cdot }_{p} \) has ultra-metric property,
the set of all discs of \( \QQ _{p} \) is countable and, if it is
equipped with the natural order given by inclusion, it has the structure
of an oriented tree: each disc of radius \( p^{k} \) contains exactly
\( p \) discs of radius \( p^{k-1} \), its sons, and the disc \( D(u,p^{k}) \)
of center \( u \) and radius \( p^{k} \) has the disc \( D(u,p^{k+1}) \)
as father. 

\vspace{0.3cm}
{\centering \resizebox*{0.8\textwidth}{!}{\includegraphics{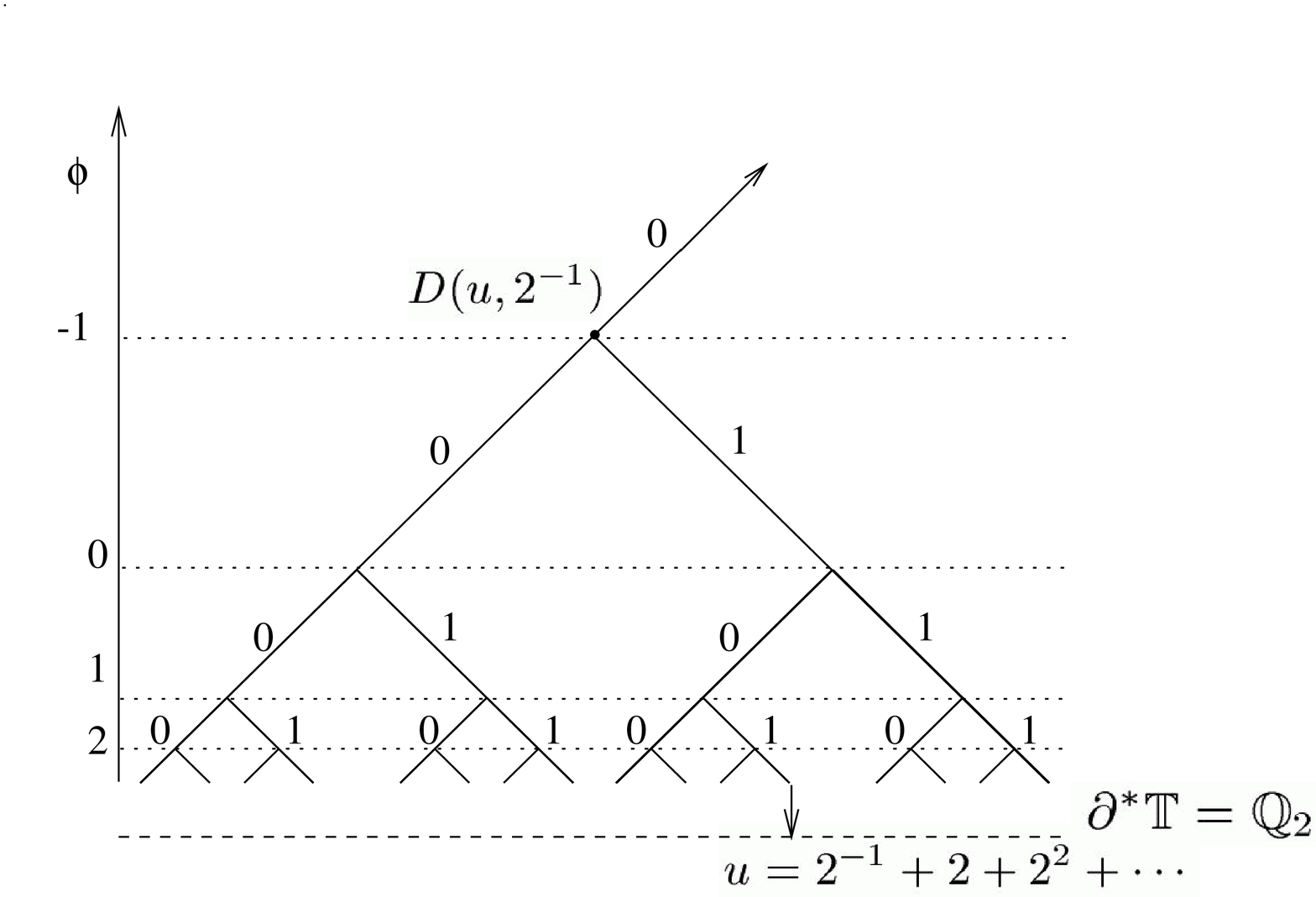}} \par}
\vspace{0.3cm}

If \( D(0,1) \) is the origin of the tree, then the Busemann function
is \[
\phi (D(u,p^{k}))=-k.\]
 One has a one to one mapping (that is in fact an isometry) between
\( \bT  \) and \( \Qp  \), associating to decreasing sequence \( \left\{ D(u,q^{-k})\right\} _{k\in \NN } \)
with the element \( u=\bigcap _{k\in \NN }D(u,q^{-k}) \) of \( \Qp  \). 
\vspace{0.3cm}

Like in the real case, the \( p \)-adic affine group, \( \affO{\Qp } \),
is the set of the mappings of the form\[
g=(t,a):\, u\mapsto au+t\quad \, \textrm{with }a\in \Qp ^{*}\textrm{ and }t\in \Qp .\]
and it can be realized as the group of matrices \[
\affO{\Qp }=\left\{ \left[ \begin{array}{cc}
a & t\\
0 & 1
\end{array}\right] :\, a\in \Qp ^{*}\textrm{ and }t\in \Qp \right\} .\]
 As affinities send discs on discs, and respect the inclusion order,
they constitute a subgroup of \( \affT  \), that is closed and non-exceptional. 

Since every affinity \( (t,a) \) is the composition of a translation
\( u\mapsto u+t \) and of a roto-homothety \( u\mapsto au \), it
is natural to see \( \affO{\Qp } \) as a semi-direct product \( \Qp \rtimes \Qp ^{*} \).
However this decomposition does not coincide with the one we have
earlier introduced as semi-direct product of \( \horO{\Qp } \) and
of \( \ZZ  \). In fact as\[
\phi ((t,a))=v_{p}(a),\]
 the horocyclic group is \[
\horO{\Qp }=\left\{ \left[ \begin{array}{cc}
a & t\\
0 & 1
\end{array}\right] \in \affO{\Qp }:\, \nor{a}=1\right\} =\Qp \rtimes \ZZ _{p}\]
where \( \ZZ _{p} \) is the ring of \( p \)-adic integers, that
coincide with the group of rotation of center \( 0 \), while \( \Qp  \)
can be identified with the group of translation.

The \( p \)-adic affine group has much more similarities with the
real case then a generic subgroup of \( \affT  \). First of all in
the \( p \)-adic setting the tree has an algebraic structure that
is much more stiff that the one of a simple graph: every affinity
is characterized by two parameters (\( t \) and \( a \)), and it
is therefore uniquely determined when one knows how it acts on two
points of \( \Qp =\bT  \). Secondly, but not less important, the
group of roto-homotheties of center \( 0 \), i.e. its stabilizer\[
K_{0}(\Qp )\times \ZZ =\Qp ^{*}\]
 is Abelian. Finally the bottom boundary \( \bT  \) is identified
with \( \Qp  \) and it has then the structure of a Abelian group.
These properties permit to obtain stronger results for random walks
on the \( p \)-adic affine group than in the general setting.

One of the interests in the \( p \)-adic affine group is linked with
the study of random walk on the group of affine transformations whose
coefficients can take only rational values, \( \affO{\QQ } \). This
group may be naturally be regarded as a dense subgroup of the real
or of the \( p \)-adic affine group, according to the metric one
considers. As it has already be pointed out in previous works (see
Kaimanovich \cite{Ka91}) from a measure theoretic point of view the
behavior of random walks on \( \affO{\QQ } \) is not necessarily
related to the Euclidean metric and a complete understanding may be
obtained by a simultaneous immersion in \( \affO{\RR } \) and in
the \( \affO{\QQ _{p}} \).

\subsubsection{Lamplighter group}

Another algebraic structure on the tree, different from the \( p \)-adic
one, but that guarantees the same regularity properties, is given
identifying the tree with sequences of integer numbers modulo \( q \),
i.e \( \ZZ /q\ZZ  \), and considering the action of the Lamplighter
group, i.e. of the wreath product \( (\ZZ /q\ZZ )\wr \ZZ  \). More
precisely let \[
\sZ _{q}=\left\{ \sigma :\ZZ \ver \ZZ /q\ZZ \, |\, \sigma \, \textrm{has finite support}\right\} \]
  and, for every \( k\in \ZZ  \), consider the equivalence relations
such that if \( \sigma _{k} \) is the class of \( \sigma \in \sZ _{q} \)
then \( \sigma _{k}=\tau _{k} \) if and only if \( \sigma (n)=\tau (n) \)
for all \( n\leq k \). One can then identify the oriented tree of
degree \( q+1 \) with the set \( \left\{ \sigma _{k}\, |\, \sigma \in \sZ _{q}\textrm{ and }k\in \ZZ \right\}  \)
in such a way that \( \sigma _{k} \) is the father of \( \sigma _{k+1} \).
The Lamplighter, that can also be seen as the semi-direct product
\( \sZ _{q}\rtimes \ZZ  \), acts on the tree the usual sum of \( \sZ _{q} \)
and by shift and it is then a non-exceptional ( but non closed) subgroup
of \( \affT  \).

\subsection{Renewal on the tree}

As we assumed that the group, \( \Gamma  \), generated by the support
of \( \mu  \) is non-exceptional, thus non-unimodular, a fundamental
result (cf. \cite{GKR}) ensures that the random walks are always
transient, i.e. almost surely they visit every compact set only a
finite number of times. Their \emph{potential measure,} i.e. \[
U(A)=\sum ^{\infty }_{n=0}\mu ^{(n)}(A)=\esp {\sum ^{\infty }_{n=0}\Ind{[L_{n}\in A]}}=\esp {\sum ^{\infty }_{n=0}\Ind{[R_{n}\in A]}},\]
is then a Radon measure on \( \Gamma  \). The \emph{(right) potential
kernel} is the family of measures\[
g*U(A)=\int _{\Gamma }\Ind{A}(gx)\, U(dx)=\esp {\sum ^{\infty }_{n=0}\Ind{[gR_{n}\in A]}}\]
where \( g\in \Gamma  \). That is the expected number of visits in
the set \( A\subset \Gamma  \) for the right random walk starting
in \( g \). By the maximum principle this family is bounded for every
compact set \( A \) when \( g \) varies in \( \Gamma  \) and thus
is vaguely relatively compact. Its limit measures when \( g \) goes
to infinity are Radon measures on the group \( \Gamma  \) that are
right \( \mu  \)-excessive, that is satisfies to the inequality \( \nu *\mu \leq \nu  \)
( they will turn out to be \( \mu  \)-invariant). 

The goal of this paper is to study these measures, to describe their
properties and to determine the directions of convergence of the potential
kernel. In section 3, we will state (and prove) in detail our main
results, that may be resumed in the following

\begin{thm*}
Suppose that the measure \( \mu  \) is spread out and that satisfies
to suitable moment conditions then the potential kernel \( g*U \)
can be continuously extended to \( \partial \Gamma  \) .

Furthermore\[
\lim _{g\ver \omega }g*U=0\]
and for all \( \beta \in \bT  \), if \( s\in \Gamma  \) is a refernce
homothety of center \( \alpha  \) and \( b\in \Gamma  \) is such
that \( b\alpha =\beta  \) \[
\lim _{g\ver \alpha }g*U=\nu _{\alpha }=\left\{ \begin{array}{l}
0\qquad \qquad \textrm{ if }\mu (\phi )>0\\
b*m_{\left\langle s\right\rangle }*\widehat{\overline{m}}\qquad \textrm{if }\mu (\phi )\leq 0
\end{array}\right. \]
where \( m_{\left\langle s\right\rangle } \) is the counting measure
on the sub-group \( \left\langle s\right\rangle =\ZZ  \) and \( \overline{m} \)
is the unique Radon \( \widehat{\mu } \)-invariant measure on \( \horg  \)
and invariant by right action of \( K_{\alpha } \) (we denote by
\( \widehat{\nu } \) the image of the measure \( \nu  \) by the
inversion on the group).
\end{thm*}
We observe the measure \( \overline{m} \) is finite if and only if
\( \mu (\phi )<0 \); in this case the total mass is \( -\frac{1}{\mu (\phi )} \).

The hypothesis that \( \mu  \) is spread out is not necessary when
\( \Gamma  \) is a sub-group of \( \affT  \) that acts in sufficiently
regular way and, in particular, when the random walk is supported
by \( \affO{\QQ _{p}} \) or by the Lamplighter group and we look
at the limits towards a point of \( \bT  \).

The characterization of the limits of potential kernel is the starting
point for a more detailed study for random walks on \( \affT  \).
Using the identification, due to Woess \cite{Wo95}, of the Martin
boundary of a random walk on a transitive subgroup of the tree isometries
with the tree boundary, we proved in \cite{tesi} that, if the measure
\( \mu  \) has a compact support and continuous density, it is possible
to give an integral representation of \( \mu  \)-invariant measures
by mean of the measures \( \nu _{\alpha } \). For instance if \( \mu (\phi )=0 \),
then every \( \mu  \)-invariant measure \( \nu  \) on \( \Gamma  \)
can be written in a unique way as\[
\nu =c_{\nu }m_{\Gamma }^{r}+\int _{\bT }\nu _{\alpha }\, \eta _{\nu }(d\alpha )\]
where \( m^{r}_{\Gamma } \) is the right Haar measure of \( \Gamma  \).
By this integral representation, in its turn, one can prove the uniqueness
of the measure \( \mu  \)-invariant on the right and the left, and,
by Guivarc'h's quotient theorems, a local limit theorem, as in the
real case (cf.\cite{LPP}). Under suitable moment and regularity conditions,
when \( \mu (\phi )=0 \), one can prove that \[
\lim _{n\ver \infty }c(n)n^{3/2}\mu ^{(n)}=\overline{m}'*m_{\left\langle s\right\rangle }*\widehat{\overline{m}}\]
where \( \overline{m}' \) is the unique \( \mu  \)-invariant measure
on \( \horg  \) and invariant by right action of \( K_{\alpha } \),
and \( c(n) \) is a sequence uniformly bounded away from \( 0 \)
and \( +\infty  \), that can be proved to be constant when \( \Gamma =\affO{\Qp } \).

We would like to conclude this section giving the guiding line of
our study. 

One of our main tools is a renewal equation, that says there exists
a probability measures \( \rho  \) such that \[
\rho *U=\nu _{\alpha }\quad \textrm{on }\left\{ g\in \Gamma \, |\, \phi (g)\leq 0\right\} ,\]
namely when the starting point of the right random is distributed
as \( \rho  \) its potential measure is given by the limit measure,
at least on half of the group. This equality is proved, after some
preliminary results, at the end of the next section.

The second main step is to determine some fundamental invariance properties
for the accumulation points of the potential kernel. In particular
we need to show that these accumulation points are left invariant:
\[
\lim _{n\ver \infty }g_{n}*U=\lim _{n\ver \infty }g_{n}g*U\qquad \forall g\in \Gamma \]
 whenever the limit exists. This is a necessary condition to show
that the potential kernel can be continuously extended to the tree
boundary, because it easily checked that whenever \( g_{n}\ver \alpha  \)
then also \( g_{n}g\ver \alpha  \). It is at this stage that there
are the main differences from the real case and where the hypothesis
that \( \mu  \) is spread out arise. In section 3 we show this and
other regularity properties and we give the proves of our main results.

\section{Preliminary results and a renewal equality}

In this section we are going to give some preliminary results. In
the first sub-section we describe the convergence of the right random
walk to the boundary improving a result originally due to Cartwright,
Kaimanovich and Woess \cite{CKW}. Next we analyze the action of the
random affinities on the bottom boundary of the tree, \( \bT  \).
Finally we provide a renewal equality, for both the action on the
boundary and the random walk on the group, that will serve as one
of our main tools. For this last result we use similar methods as
those that have been used for the study of random walks of the real
affine group by Babillot, Bougerol and Elie in \cite{BBE}.

\subsection{Convergence of the random walk to the boundary}

Since the closed subgroup generated by the support of the law \( \mu  \)
is non-exceptional, the random walks \( L_{n} \) and \( R_{n} \)
are transient so that the accumulation points of their trajectories
lie on the boundary of the tree. Namely the right random walk converges
to a random variable on the boundary \( \partial \TT  \). 

\begin{thm}
\label{teofond} Suppose that \( \esp {\nor{\phi (X_{1})}}<\infty  \). 
\begin{enumerate}
\item If \( \mu (\phi )<0 \) then \( R_{n}\rightarrow \omega  \) almost
surely. 
\item If \( \mu (\phi )>0 \) and \( \esp {\nor{X_{1}}}<\infty  \) then
\( R_{n}\rightarrow \xi _{\infty } \) almost surely where \( \xi _{\infty } \)
is a random element in \( \bT  \). The law \( m \) of \( \xi _{\infty } \)
is supported by \( \bG  \) and carries no point mass.
\item If \( \mu (\phi )=0 \) and \( \esp {\nor{X_{1}}}<\infty  \) then
\( R_{n}\rightarrow \omega  \) almost surely.
\end{enumerate}
\end{thm}
\begin{proof}
Results (1) and (2) are in \cite{CKW}, who proved (3) only under
an exponential moment condition. Our proof of (3) uses a method developed
in \cite{Br03}. 

For all \( x\in \TT  \), consider the cone \[
C_{x}=\left\{ y\in \TT :x\succeq y\right\} .\]
We show that for every \( x\in \TT  \)\[
\pr {R_{n}o\in C_{x}\, \textrm{infinitely often }}=0.\]
Let \( m^{r}_{\Gamma } \) be the right Haar measure of \( \Gamma  \).
First, we prove that for all \( x\in \TT  \) and \( m^{r}_{\Gamma } \)-almost
all \( g\in \Gamma  \)\[
\pr {gR_{n+1}o\in C_{x},\, gR_{n}o\not \in C_{x}\, \textrm{ infinitely often }}=0.\]
 Using the Borel-Cantelli Lemma, it is sufficient to show that \[
\sum ^{\infty }_{n=0}\pr {gR_{n+1}o\in C_{x},\, gR_{n}o\not \in C_{x}}=g*U(\psi )<+\infty \]
where \( \psi (g)=\pr {gX_{1}o\in C_{x},\, go\not \in C_{x}} \),
and using Lemma 2.2 in \cite{Br03} we just need to show that \( \psi  \)
is \( m^{r}_{\Gamma } \)-integrable.

Observe that\begin{eqnarray*}
\int _{\Gamma }\psi (g)m_{\Gamma }^{r}(dg) & = & \esp {\int _{\Gamma }\Ind{[gX_{1}o\in C_{x},\, go\not \in C_{x}]}m_{\Gamma }^{r}(dg)}\\
 & = & \esp {\int _{\Gamma }\Ind{[go\in C_{x},gX^{-1}_{1}o\not \in C_{x}]}m_{\Gamma }^{r}(dg)}
\end{eqnarray*}
and that\begin{eqnarray*}
\left\{ g\in \Gamma :\, go\in C_{x},gX^{-1}_{1}o\not \in C_{x}\right\}  & = & \left\{ g\in \Gamma :\, go\in C_{x},g(o\cf X^{-1}_{1}o)\not \in C_{x}\right\} \\
 & = & \bigcup _{y\in S}V(y\ver o)
\end{eqnarray*}
where \( V(y\ver o)=\left\{ g\in \Gamma :\, gy=o\right\}  \) and
\( S \) is the geodesic segment that joins \( o \) to the vertex
immediately before \( o\cf X^{-1}_{1}o \) (or the empty set if \( o\cf X^{-1}_{1}o=o \)).
Observe that for any affinity \( \gamma  \) such that \( \gamma y=o \)
one has \( V(y\ver o)=V(o\ver o)\gamma  \), thus \( m^{r}_{\Gamma }(V(o\ver o))=m^{r}_{\Gamma }(V(y\ver o)) \)
for all \( y \). Since the segment \( S \) contains exactly \( -\phi (o\cf X_{1}^{-1}o) \)
vertices, we have\[
\int _{\Gamma }\Ind{[go\in C_{x},gX^{-1}_{1}o\not \in C_{x}]}m_{\Gamma }^{r}(dg)\leq -\phi (o\cf X_{1}^{-1}o)m^{r}_{\Gamma }(V(o\ver o))\]
 so that\[
\int _{\Gamma }\psi (g)m_{\Gamma }^{r}(dg)\leq \esp {-\phi (o\cf X_{1}^{-1}o)}m^{r}_{\Gamma }(V(o\ver o))\leq \esp {\left| X_{1}\right| }m^{r}_{\Gamma }(V(o\ver o))<+\infty .\]
We proved that for almost all \( g \), almost surely, \( gR_{n}o \)
cannot pass from \( C_{x}^{c} \) to \( C_{x} \) but a finite number
of times.

As the set \( V(x\ver x) \) is open and has then strictly positive
Haar measure, there exists \( g\in V(x\ver x) \) such that the event
\begin{eqnarray*}
\left[ gR_{n+1}o\in C_{x},\, gR_{n}o\not \in C_{x}\right]  & = & \left[ R_{n+1}o\in g^{-1}C_{x},\, R_{n}o\not \in g^{-1}C_{x}\right] \\
 & = & \left[ R_{n+1}o\in C_{x},\, R_{n}o\not \in C_{x}\right] .
\end{eqnarray*}
 takes place only a finite number.

On the other hand the vertex set of the tree is countable, whence
almost surely\[
\forall x\in \TT :\, \, R_{n+1}o\in C_{x},\, R_{n}o\not \in C_{x}\textrm{ }\quad \textrm{a finite number of times}.\]

If \( \mu (\phi )=0 \), the random walk \( \phi (R_{n}) \) is recurrent
on \( \ZZ  \), and visits the interval \( ]-\infty ,\phi (x)] \)
infinitely often. In particular \( R_{n}o\in C^{c}_{x} \) infinitely
often. As \( R_{n} \) cannot go back to \( C_{x} \) but a finite
number of times, we conclude that it is in \( C^{c}_{x} \) for all
sufficiently large \( n \) . 
\end{proof}

\subsection{Action on \protect\( \bT \protect \)}

As we have noted, the affine group has a natural action on the boundary
of the tree that may also be considered as a boundary for the group
itself; the behavior of the Markov chain on \( \bT  \) provides then
useful information for the study of the random walk on the group itself.

The analogue of this Markov chain in the setting of real affine group
is the process induced on \( \RR  \) by the natural action. We may
note that while in the real case the joint behavior of this chain
and of the homothetic component describes completely the random walk,
this does not hold on the tree because the process on the group is
more complex. In fact the \( \affT  \) is the semi-direct product
of \( \ZZ  \) and of the horocyclic group \( \hor  \), that is {}``bigger''
than the boundary of the tree. 

To understand the behavior of the random walk on the group we also
need to analyze the process induced on the horocyclic group by the
action of \( \affT  \). However the techniques we are going to use
cannot be applied directly in this setting, mainly because the action
of \( \ZZ  \) on the horocyclic group is not sufficiently contractive.
However at the end of this section we will be able to deduce, as a
corollary, some results in this context as well.

Let \( \Upsilon _{0} \) be a random variable defined on the probability
space \( (\Omega ,\PP ) \), with value in \( \bT  \), independent
from the increments \( \{X_{n}\}_{n} \) of the random walk. The Markov
chain induced on \( \bT  \) is the process 

\[
\Upsilon _{k}=L_{k}\Upsilon _{0}=X_{k}\cdots X_{1}\Upsilon _{0}.\]
Its transition kernel is \[
Pf(\upsilon )=\esp {f(X_{1}\cdot \upsilon )}=\int _{\affT }f(g\upsilon )\mu (dg)=\mu \stackrel{\cdot }{*}\upsilon (f).\]
(Here the symbol \( \stackrel{\cdot }{*} \) denotes the convolution
of a measure on the group \( \Gamma  \) and a measure on a \( \Gamma  \)-space). 

The proofs of some results of this section are formally very similar
to the analogues in the real case. We have translated them in this
setting for readers convienence and they can be found in the appendix.

The behavior of this chain is directly related to the random walk
\( S_{n}=\phi (L_{n}) \) on \( \ZZ  \), that contains the information
on how the random affinities contract or dilate the boundary. 

When \( \mu (\phi )=\esp {\phi (X_{1})}\not =0 \), one may directly
obtain properties of transience or recurrence of the induced Markov
chain from what is known for the process on the group. 

\begin{prop}
\label{prop-munon0} Suppose that \( \esp {\nor{X_{1}}}<\infty  \).
\begin{enumerate}
\item If \( \mu (\phi )<0 \) then for all \( \upsilon \in \bT  \) \[
\lim _{n\rightarrow \infty }L_{n}\upsilon =\omega \qquad \textrm{almost surely}\]
 and the chain \( \left\{ \Upsilon _{n}\right\} _{n} \) is transient. 
\item If \( \mu (\phi )>0 \), the law \( m \) of the random variable \( \ds \xi _{\infty }=\lim _{n\rightarrow \infty }R_{n} \)
is the unique probability measure for the Markov chain on \( \bT  \),
thus \( \left\{ \Upsilon _{n}\right\} _{n} \) is positive recurrent.
Furthermore for all \( \upsilon  \) in \( \bT  \) \( \left\{ L_{n}\upsilon \right\} _{n} \)
visits infinitely often every open set of \( \bT  \) of non null
\( m \)-measure, almost surely.
\end{enumerate}
\end{prop}
\begin{proof}
see Appendix
\end{proof}
A classical technique to deal with the centered case, \( \mu (\phi )=0 \),
is to extract from the original chain a sub-chain with positive drift,
to which one can apply the previous results. We consider the sequence
of ladder stopping times \( l_{k} \) giving the times when the random
walk \( S_{n}=\phi (L_{n}) \) on \( \ZZ  \) reaches a new maximum\[
l_{k}=\min \{n>l_{k-1}:\, \, S_{n}>S_{l_{k-1}}\}\quad \textrm{ and }\quad l_{0}=0.\]
 The process obtained regarding the left random walk at these times
\[
L_{l_{k}}=\left( X_{l_{k}}\cdots X_{l_{k-1}+1}\right) L_{l_{k-1}}\]
is still a left random walk, whose law, that is the law of \( L_{l_{1}} \),
will be denoted by \( \mu _{l} \). For convenience we also set \( l_{1}=l \).

The drift of this random walk is clearly positive (eventually infinite)
\[
\esp {\phi (L_{l})}>0,\]
 but to apply the previous results one has to make sure that its law
is integrable and for this we need a moment of order \( 2+\varepsilon  \)
.

\begin{lem}
\label{lem-moment.marche.descend} The closed group generated by the
support of the law of \( L_{l_{1}} \) coincide with the group generated
by the support of \( \mu  \) and if one assume that \[
\esp {\phi (X_{1})^{2}+\nor{b(X_{1})}^{2+\varepsilon }}<+\infty \]
 then \[
\esp {\nor{L_{l}}}<+\infty .\]

\end{lem}
\begin{proof}
The result is proved in \cite{CKW} (Proposition 4), with the slight
difference that there it is formulated for right random walk and decreasing
ladder times. 
\end{proof}
Under this stronger moment hypothesis, we can apply to the random
walk \( L_{l_{n}} \) the results of the previous proposition, in
particular to ensure that the Markov chain \( \Upsilon _{l_{n}} \)
has a unique invariant probability measure that is denoted \( m_{l} \).
We get then the following results:

\begin{prop}
\label{prop-Ric.mu0}Suppose that \( \mu  \) has a moment of order
\( 2+\varepsilon  \) and that \( \mu (\phi )=0 \). Then the chain
\( \Upsilon _{n} \) is recurrent, in the sense that almost surely
for all \( \upsilon  \) in \( \bT  \) the chain \( L_{n}\upsilon  \)
visits infinitely often every open set of non null \( m_{l} \)-measure.
Furthermore there exists a unique \( \mu  \)-invariant Radon measure
on \( \bT  \). 
\end{prop}
\begin{proof}
Applying to \( L_{l_{n}}\upsilon  \) the results of Proposition \ref{prop-munon0},
we obtain that \( L_{l_{n}}\upsilon  \) visits infinitely often every
open set of non null \( m_{l} \)-measure, whence, a fortiori, the
same holds for \( L_{n}\upsilon  \). 

Thus for every \( \upsilon \in \bT  \) and any non-negative continuous
function \( f \) on \( \bT  \) such that \( m_{l}(f)\neq 0 \)\[
\sum ^{+\infty }_{n=0}P^{n}f(\upsilon )=\esp {f(L_{n}\cdot \upsilon )}=+\infty \]
 and the chain is topologically conservative. As \( P \) is a Feller
operator, by \cite{Lin70} (Theorem 5.1), the chain \( \Upsilon _{n} \)
has an invariant Radon measure.

In the next lemma we shall prove what is known in as local contraction
property. Using this result and Chacon-Ornstein Theorem one obtains
uniqueness of the invariant measure along the same lines in \cite{Br03},
Theorem 3.
\end{proof}
In the centered case the random walk \( L_{n} \) has neither a contracting
or dilating action, as the distance between two trajectories \[
\Theta (L_{n}\upsilon ,L_{n}\varsigma )=q^{-\phi (L_{n})}\Theta (\upsilon ,\varsigma )\]
does not converge to zero as in the case of positive drift , nor to
\( +\infty  \) as in the case of negative drift, but oscillates between
these two extremes. However, if we do not look globally at this process
but only through a compact window, we can recover a stability property:

\begin{lem}
If the measure \( \mu  \) has a first moment and if \( \mu (\phi )=0 \),
for all compact set \( K \) in \( \bT  \) and for every pair of
ends \( \upsilon ,\varsigma \in \bT  \) , almost surely\begin{equation}
\label{eq-cont-loc}
\lim _{n\rightarrow \infty }\Theta (L_{n}\upsilon ,L_{n}\varsigma )\, \Ind{K}(L_{n}\upsilon )=\lim _{n\rightarrow \infty }\Theta (\upsilon ,\varsigma )\, q^{-\phi (L_{n})}\Ind{K}(L_{n}\upsilon )=0.
\end{equation}

\end{lem}
\begin{proof}
The local contraction (\ref{eq-cont-loc}) is a direct consequence
of the fact that, in the centered case, the right random on \( \affT  \)
converges to the mythical ancestor \( \omega  \). 

Suppose that (\ref{eq-cont-loc}) does not hold. Then there exists
a cone with vertex \( y\in \TT  \) \[
C_{y}=\left\{ x\in \TT \cup \bT :x\succeq y\right\} \]
and an integer \( M \) such that, with probability 1,\[
L_{n}\upsilon \in C_{y}\textrm{ and }\phi (L_{n})<M\textrm{ }\quad \textrm{for infinitely many }n.\]
Then for any \( x \) in the geodesic \( \overline{\upsilon \omega } \)
such that \( \phi (x)<\phi (y)-M \) one has \begin{eqnarray*}
1=\pr {L_{n}x\in \overline{y\omega }\textrm{ infinitely often}} & = & \pr {x\in \overline{\left( L_{n}\inv y\right) \omega }\textrm{ infinitely often}}\\
 & \leq  & \pr {L_{n}\inv y\in C_{x}\textrm{ infinitely often}}.
\end{eqnarray*}
On the other hand \[
L_{n}^{-1}=X_{n}^{-1}\cdots X_{1}^{-1}=\hat{R}_{n}\]
 is a right random walk with first moment and null drift. Thus we
have obtained a contradiction, because we know by Theorem \ref{teofond}
that \( L_{n}\inv =\hat{R}_{n}\rightarrow \omega  \) almost surely.
\end{proof}
As in the real case, it is possible to construct the unique invariant
Radon measure, \( m \), using the invariant probability measure,
\( m_{l} \), of the contracting sub-chain, in the following way:
\begin{equation}
\label{eq-misinv}
m(f)=\frac{1}{\esp {S_{l}}}\int _{\bT }\esp {\sum _{k=0}^{l-1}f(L_{k}\cdot \upsilon )}m_{l}(d\upsilon )
\end{equation}
Using the strong Markov property one can see that \( m \) is \( \mu  \)-invariant.
Observe that the stopping time \( l \) is not integrable (cf. \cite{Spi}),
so that the measure \( m \) does not have finite mass. The fact that
it is finite on compact sets is not evident and it will proved in
the next sub-section in Corollary \ref{cor-mis.Rad}. 

The measure (\ref{eq-misinv}) can also be defined if the step law
of the random walk has first moment and the drift is positive. In
this case the time \( l \) is integrable and \( m \) has total mass
equal to\[
m(\bT )=\frac{\esp {l}}{\esp {S_{l}}}=\frac{\esp {l}}{\esp {l}\esp {S_{1}}}=\frac{1}{\mu (\phi )}.\]
by Wald's equality (cf. \cite{Doo}, page 350). This measure is then
just a normalization of the unique invariant probability measure.

\subsection{A renewal equality}

In order to obtain the announced renewal equality we consider the
joint action of random affinities on both the bottom boundary and
the integers. In other words the result concerns the Markov chain
on \( \bT \times \ZZ  \) whose transition kernel is \[
\tilde{P}f(\upsilon ,z)=\mu \stackrel{\cdot }{*}(\upsilon ,z)(f)=\esp {f\left( X_{1}\cdot (\upsilon ,z)\right) }=\esp {f(X_{1}\upsilon ,\phi (X_{1})+z)}\]
and whose potential kernel is\[
\sum _{n=0}^{\infty }\tilde{P}^{n}f(\upsilon ,z)=U\stackrel{\cdot }{*}(\upsilon ,z)(f)\]
 where \( U=\sum _{n=0}^{\infty }\mu ^{(n)} \) is the potential measure
of the random walk on the group. 

\begin{prop}
\label{prop-Pot.bT=3Dmxl}Assume that \( \mu (\phi )>0 \) and a moment
of first order or that \( \mu (\phi )=0 \) and a moment of order
\( 2+\varepsilon  \). Let \( m \) be the measure defined in (\ref{eq-misinv}).
There exists a probability measure \( p \) on \( \bT \times \ZZ  \)
such that for every non-negative function \( f \) with support in
\( \bT \times \ZZ _{+} \) \begin{equation}
\label{eq-mipot}
U\stackrel{\cdot }{*}p(f)=m\times m_{\ZZ }(f)
\end{equation}
 where \( m_{\ZZ } \) is the counting measure on \( \ZZ  \).
\end{prop}
\begin{proof}
see Appendix 
\end{proof}
As announced, a direct consequence of this result is the following 

\begin{cor}
\label{cor-mis.Rad}Under the hypothesis of the previous proposition,
the measure \( m \) defined in (\ref{eq-misinv}) is a Radon measure
and coincides with the unique \( \mu  \)-invariant Radon measure
on \( \bT  \).
\end{cor}
\begin{proof}
As the random walk on the affine group is transient, its left potential
kernel \( U*g(f) \) is bounded for every bounded function \( f \)
with compact support. As the group \( \Gamma =\horg \rtimes \ZZ  \)
is a compact extension of \( \bT \times \ZZ  \) , it easy to see
that also the potential kernel of the chain induced on \( \bT \times \ZZ  \)
that is \[
U\stackrel{\cdot }{*}x(f)=\int _{\Gamma }f(g\cdot x)U(dg)\]
is bounded in \( x\in \bT \times \ZZ  \) for any bounded compactly
supported function \( f \).

Let \( K \) be a compact set of \( \bT  \), then, as \( p \) is
a probability, one has \[
m(K)=m\times m_{\ZZ }(\Ind{K}\times \Ind{\left\{ 0\right\} })=\int _{\bT \times \ZZ }U\stackrel{\cdot }{*}x(\Ind{K}\times \Ind{\left\{ 0\right\} })p(dx)<+\infty \]

\end{proof}

\subsection{Action on \protect\( \horg \protect \) and a renewal equality on
the group}

To understand the random walk on the group \( \Gamma  \), we are
also interested on the process induced on the horocyclic group by
the action defining the semi-direct product \( \Gamma =\horg \rtimes _{s}\ZZ  \),
that is \[
g\cdot b:=b(g)s^{\phi (g)}bs^{-\phi (g)}=gbs^{-\phi (g)}.\]
 Even if the techniques that we have used are not completely adapted
to this setting, we can easily deduce some useful results.

As we have observed, the bottom boundary of the tree is homeomorphic
to the quotient of the horocyclic group by the compact stabilizer
of an end \( \alpha  \), that is \( \bG =\horg /K_{\alpha } \).
It is then possible to extend any measure, \( m \), on \( \bG  \)
to a measure, \( \overline{m} \), on \( \horg  \) by setting \begin{equation}
\label{eq-est.mis}
\overline{m}(f)=\int _{\bG \times K_{\alpha }}f(xk)m_{K_{\alpha }}(dk)m(dx)
\end{equation}
 where \( m_{K_{\alpha }} \) is the Haar measure of \( K_{\alpha } \).
By construction \( \overline{m} \) is right invariant by the right
action of the rotations \( K_{\alpha } \) and, above all, it is well
adapted to the actions of the group \( \Gamma  \) on the boundary
and on the horocyclic group, respectively. In fact for every affinity
\( g \), one has\[
g\stackrel{\cdot }{*}\overline{m}=\overline{g\stackrel{\cdot }{*}m}.\]
Thus, whenever one has an invariant measure on the boundary, its possible
to extend it to an invariant measure on the horocyclic group (right
invariant by action of \( K_{\alpha } \)). It is also possible to
extend the renewal equality of Proposition \ref{prop-Pot.bT=3Dmxl}
to \( \horg \times \ZZ  \), in order to get a renewal equality on
the whole group. To resume, we obtain the following

\begin{cor}
\label{cor-eq.ren.AffT} Assume that \( \mu (\phi )>0 \) and a moment
of first order or that \( \mu (\phi )=0 \) and a moment of order
\( 2+\varepsilon  \). Then there exists a \( \mu  \)-invariant Radon
measure, \( \overline{m} \), on \( \horg  \) obtained by extension
of the invariant measure on \( \bG  \). This measure has finite mass
if and only if \( \mu (\phi )>0 \). 

If \( \overline{m} \) is normalized in such a way that it is the
extension of the measure defined (\ref{eq-misinv}), than there exists
a probability \( \overline{p} \) on \( \Gamma  \) such that \begin{equation}
\label{eq-ren.AffT}
U*\overline{p}=\overline{m}\times m_{\ZZ }=\overline{m}*m_{\left\langle s\right\rangle }\quad \textrm{on }\left\{ g\in \Gamma \, |\, \phi (g)\geq 0\right\} ,
\end{equation}
where \( m_{\left\langle s\right\rangle } \) is the counting measure
on the group \( \left\langle s\right\rangle =\ZZ  \).
\end{cor}
\begin{proof}
Just observe that equation (\ref{eq-ren.AffT}) can be obtained proceeding
as in Proposition \ref{prop-Pot.bT=3Dmxl}, where we never really
used the structure of the chain on the boundary, but just the fact
that there exists a probability measure for the contracting sub-chain.
Otherwise, we can obtain (\ref{eq-ren.AffT}) by extension from its
equivalent in the boundary case by setting\[
\overline{p}(f)=\int _{\bT \times K_{\alpha }\times \ZZ }f(xks^{h})m_{K_{\alpha }}(dk)p(dx\, dh).\]

\end{proof}
The question that remains still open is whether the \( \mu  \)-invariant
measure defined by extension is the unique invariant measure, or,
in other words, if all invariant measures are invariant by \( K_{\alpha } \).
In the case of \( p \)-adic affine group, it seems likely to have
a positive answer, that could be obtained regarding the random walk
on the group of rotations obtained by projection. In the general case
there is no such canonical projection on \( K_{\alpha } \) and the
solution do not seem to be so evident.

\begin{rem*}
As we did not impose continuity conditions on the law of the random
walk, we can apply these results to products of random affine transformations
with rational coefficients, considered as random walks on the group
\( \mathrm{Aff}(\QQ _{p}) \). One obtains the following 
\end{rem*}
\begin{cor}
\label{cor-lim.pot.QQ}Let \( \mu  \) be a measure on \( \affO{\QQ }=\QQ \rtimes \QQ ^{*} \).
Let \( t \) and \( a \) be the projections of \( \affO{\QQ } \)
on \( \QQ  \) and \( \QQ ^{*} \) respectively. Assume the measure
\( \mu  \) is irreducible on \( \affO{\QQ _{p}} \) i.e. \[
\pr {\nor{a(X_{1})}_{p}=1}<1\quad \textrm{and}\quad \forall y\in \QQ :\, \, \, \, \: \pr {a(X_{1})y+t(X_{1})=y}<1\]
 and that moment conditions of the preceding theorems are satisfied. 

Then if \( \esp {\log \nor{a(X_{1})}_{p}}\leq 0 \) there exists a
unique \( \mu  \)-invariant Radon measure on \( \Qp  \), that has
finite mass if and only if \( \esp {\log \nor{a(X_{1})}_{p}}<0 \).
\end{cor}
It is easily checked the measure \( \overline{m}\times m_{\ZZ } \)
is left \( \mu  \)-invariant on the group. Thus for every integrable
function \( f \), the functions \[
h((t,a))=(t,a)*(\overline{m}\times m_{\ZZ })(f)\]
 are non-trivial right harmonic on \( \affO{\QQ } \), and they are
bounded if \( \esp {\log \nor{a(X_{1})}_{p}}<0 \) and \( f \) is
bounded.

The group of affine transformations with rational coefficients is
usually regarded as a dense subgroup of the group of affine transformations
with real coefficient, \( \affO{\RR } \). However, generic random
walks on \( \affO{\QQ } \) do not always behave as random walks that
are completely adapted to the topology of \( \affO{\RR } \), like
those whose law is spread out on the latter group. For instance it
is known that in this last case when \( \esp {\log \nor{a(X_{1})}}<0 \),
there is no bounded harmonic function, i.e the Poisson boundary is
trivial. On the other hand it has been shown (\cite{Ka91}) that for
random walks on the affine group of dyadic integers such that \[
\esp {\log \nor{a(X_{1})}}=-\esp {\log \nor{a(X_{1})}_{2}}<0\]
the Poisson boundary is not-trivial and coincide with \( \QQ _{2} \).
The last corollary shows that every random walks on \( \affO{\QQ } \)
whose \( p \)-adic drift is negative (for some \( p \)) has a non
trivial Poisson boundary (independently from the real drift). It seems
likely that a complete description of the Poisson boundary of the
random walk on \( \affO{\QQ } \) can be obtained immersing this group
simultaneously in the real affine group and in all the \( p \)-adic
ones.

\section{Limit measures of the potential kernel}

This section is devoted to the study of the limit measures of the
potential kernel near the boundary of the group. W.Woess \cite{Wo95}
studied the asymptotic behavior of the Martin kernel (i.e. the normalization
of the potential kernel) for a random walks, with continuous and compactly
supported density, on any closed and transitive group of isometries
of homogeneous tree and showed that the Martin boundary can be identified
with \( \partial \TT  \). For the potential kernel of random walks
on \( \affT  \) it is possible to obtain the same kind of results
of continuous extension to the geometrical boundary without asking
for a compact support and moreover to obtain the form of the limit
measures in terms of the invariant measure on the boundary and of
the counting measure on \( \ZZ  \).

In comparison with Elie's works on Lie groups, by which our study
was originally inspired, we have to face, as we announced, some new
phenomena of non-commutativity and the lack of stiffness of the tree
structure. The use of the renewal equality introduced in the last
section enables us to appreciably simplify the proofs and, mostly,
to avoid continuity hypotheses for random walks on sufficiently regular
groups, as \( \affO{\Qp } \), when one looks to the limit of the
potential kernel towards a point of \( \bT  \). 

Its also worth observing that all the conclusions of this papers can
be obtained in the same way for the real case. In particular we can
improve Elie's results assuring that even if the measure is not spread
out, the associated potential kernel on on the real affine group,
\( (t,a)*U \), converges to a limit measure whenever \( (t,a) \)
converges to \( (t_{0},0) \). 

In our study, we will need the following general results on uniform
continuity of the potential kernel (cf. \cite{El82}, Proposition
2.7 and Theorem 2.9) 

\begin{lem}
\label{lem-contuni} Let \( \{g_{n}\}_{n} \) be a sequence of elements
in \( \Gamma  \) such that \( \left\{ g_{n}*U\right\} _{n} \) vaguely
converges.
\begin{enumerate}
\item \emph{Left continuity}. There exists a subsequence \( \{g_{n_{k}}\}_{k} \)
such that, for all \( y\in \Gamma  \), the measure sequences \( \left\{ yg_{n_{k}}*U\right\} _{k} \)
vaguely converge. Furthermore, for all \( f\in C_{c}(\Gamma ) \),
the sequences \( \left\{ yg_{n_{k}}*U(f)\right\} _{k} \) converge
uniformly when \( y \) is in a compact set.
\item \emph{Right continuity}. If \( \mu  \) is spread out, there exists
a subsequence \( \{g_{n_{k}}\}_{k} \) such that, for all \( y\in \Gamma  \),
the measure sequences \( \left\{ g_{n_{k}}y*U\right\} _{k} \) vaguely
converge. Furthermore, for all \( f\in C_{c}(\Gamma ) \), the sequences
\( \left\{ g_{n_{k}}y*U(f)\right\} _{k} \) converge uniformly when
\( y \) is in a compact set.
\end{enumerate}
\end{lem}
We observe that left continuity is almost straightforward, while right
continuity is a more subtle phenomenon, that requires a stronger regularity
condition (the measure is spread out). We will often use this Lemma
to guarantee that when we perturb the sequence \( g_{n} \) (to the
right or to the left) by a sequence \( y_{n} \) that converges to
\( y \), then on a sub sequence the limit does not change if we replace
\( y_{n} \) by \( y \); for instance \begin{equation}
\label{eq-LimUnifNMart}
\lim _{k\rightarrow \infty }{g_{n_{k}}y_{n_{k}}}*U(f)=\lim _{k\rightarrow \infty }{g_{n_{k}}y}*U(f).
\end{equation}

As the structure of the limit measures in a neighborhood of \( \omega  \)
differs from the structure of the limit measure in the neighborhood
of a point in the bottom boundary \( \bG  \), we will study the two
cases in two different sub-sections.

\subsection{Limits near \protect\( \bT \protect \)}

We start by showing some invariance properties for the accumulation
point of the potential kernel. In a second step using these results
and the renewal formula of the previous section, we show that it is
possible to extend by continuity the potential kernel to \( \bG  \).

\subsubsection{Invariance properties}

The structure of the limit measures depends not only on how the affinities
act on the vertices at finite distance, but also on how they act near
the boundary of the tree. Observe that while the first of these actions
is completely adapted to the topology of pointwise convergence on
the tree), the second action is not at all related to it: one can
construct a sequence of elements of \( \affT  \) that is closer and
closer to the identity, but such that they act non trivially on a
sequence of sets sufficiently far away from the origin. For the subgroups
of \( \affT  \) that have homogeneous action on the tree (as \( \affO{\QQ _{p}} \))
this problem does not occur, but for more general subgroups to link
the topology of the group and the action near the boundary will require
the right continuity of the potential kernel (and therefore the spread
out hypothesis on \( \mu  \) ). 

The following lemmas are intended to clarify the behavior of an affinity
near the boundary, and more precisely what happens if it is conjugates
with a sequence of transformations that let one see how it acts faraway
from an end \( \alpha  \).

\begin{lem}
\label{lem-conv.alpha}Let \( \alpha \in \bT  \) be a fixed end of
the tree and let \( s\in \Gamma  \) be a reference homothety such
that \( s\alpha =\alpha  \) and \( \phi (s)=1 \). Then, for every
\( g\in \Gamma  \), the sequence \( \left\{ s^{n}gs^{-n}\right\} _{n\in \NN } \)
is relatively compact and all its accumulation points fix \( \alpha  \).
Furthermore, if \( g\in \horg  \) then the accumulation points belong
to \( K_{\alpha } \).
\end{lem}
\begin{proof}
We first observe that \( \left\{ s^{n}gs^{-n}\alpha \right\} _{n} \)
converges to the end \( \alpha  \), in fact \begin{eqnarray}
\lim _{n\toinf }\Theta (\alpha ,s^{n}gs^{-n}\alpha ) & = & \lim _{n\toinf }\Theta (s^{n}\alpha ,s^{n}g\alpha )\qquad \textrm{ because }s^{n}\alpha =\alpha \nonumber \\
 & = & \lim _{n\toinf }q^{-n}\Theta (\alpha ,g\alpha )\label{eq-ConvAlpha} \\
 & = & 0.\nonumber 
\end{eqnarray}
Since \( \phi (s^{n}gs^{-n})=\phi (g) \), the elements of the sequence
\( \left\{ s^{n}gs^{-n}\right\} _{n} \) belong, for every sufficiently
large \( n \), to the compact sets of \( \Gamma  \) of the form\[
V(\alpha _{m}\rightarrow \alpha _{m+\phi (g)})=\left\{ \gamma \in \Gamma :\, \gamma \alpha _{m}=\alpha _{m+\phi (g)}\right\} \]
for every integer \( m \), where \( \alpha _{m} \) is the vertex
of the geodesic \( \overline{\alpha \omega } \) such the \( \phi (\alpha _{m})=m \).
Thus every accumulation point of \( \left\{ s^{n}gs^{-n}\right\} _{n} \)
belongs to \[
\bigcap _{m\in \ZZ }V(\alpha _{m}\rightarrow \alpha _{m+\phi (g)})=\bigcap _{m\in \ZZ }s^{\phi (g)}V(\alpha _{m}\rightarrow \alpha _{m})=s^{\phi (g)}K_{\alpha }\]
and it fixes the end \( \alpha  \). 
\end{proof}
If \( \Gamma  \) is a subgroup that acts homogeneously on the tree,
then the sequence \( s^{n}gs^{-n} \) converges and its limit is given
by the rotation and homothetic component of center \( \alpha  \)
of \( g \). For instance if \( \Gamma =\affO{\QQ _{p}}=\QQ ^{*}_{p}\times \QQ _{p} \),
and we chose as reference homothety \( s=(p,0) \), then for \( g=(a,x) \),
the sequence \( \left\{ s^{n}gs^{-n}\right\} _{n} \) converges to
\( (a,0) \). This does not hold for general subgroup, because the
homothety \( s \) does not commute with the rotations, and this constitutes
one of the main differences with the study of the renewal on Lie groups.
L.Elie in \cite{El82} has in fact shown that any almost connected
Lie groups whose potential kernel has an infinite number of accumulation
points is the semi-direct product of a nilpotent Group (of {}``translations'')
and of the \emph{direct} product of a compact Lie group (the rotations)
and \( \RR  \) (the homotheties).

In the following lemma enables to get round the problem that for \( b\in \horg  \)
the sequence \( \left\{ s^{n}bs^{-n}\right\} _{n} \) does not converge
by giving an approximation of its accumulation points with sequences
of the form \( \left\{ s^{n_{k}}r_{k}s^{-n_{k}}\right\} _{k} \),
where \( \left\{ r_{k}\right\} _{k} \) is a sequence that converges
to a rotation \( r \). Roughly speaking \( r \) tells how \( b \)
acts far away from \( \alpha  \), but this information depends on
the subsequence \( n_{k} \) (i.e. on the escape speed from \( \alpha  \))
and, in a general case, it is not related on how \( b \) behaves
near \( \alpha  \).

\begin{lem}
Let \( \mathbf{n}=\{n_{k}\}_{k\in \NN } \) be a sequence in \( \ZZ  \)
such that \( \lim _{k\toinf }n_{k}=+\infty  \). Then for every \( b\in \horg  \)
there exists a sequence \( \left\{ r_{l}\right\} _{l} \) in \( \horg  \)
that converges to a rotation \( r \) in \( K_{\alpha } \) and a
subsequence \( \mathbf{m}=\left\{ m_{l}\right\} _{l\in \ZZ } \) of
\( \mathbf{n} \) such that\[
\lim _{l\toinf }s^{m_{l}}b^{-1}r_{l}s^{-m_{l}}=e\]

\end{lem}
\begin{proof}
By the preceding lemma one knows that there is a subsequence \( \mathbf{n}' \)
of \( \mathbf{n} \) such that \( \left\{ s^{n'_{k}}bs^{-n'_{k}}\right\} _{k} \)
converges and, therefore, is a Cauchy sequence, i.e.: \[
\lim _{k\toinf }\left( s^{n'_{k}}bs^{-n'_{k}}\right) ^{-1}s^{n'_{k+i}}bs^{-n'_{k+i}}=e\]
 uniformly in \( i\in \NN  \). Let \( \left\{ i_{k}\right\} _{k\in \NN } \)
be a sequence of natural numbers such that \( n'_{k+i_{k}}-n'_{k} \)
converges to \( +\infty  \), when \( k \) goes to \( +\infty  \).
According to Lemma \ref{lem-conv.alpha} the sequence \[
b_{k}=s^{n'_{k+i_{k}}-n'_{k}}bs^{-(n'_{k+i_{k}}-n'_{k})}\]
 is relatively compact, so that it is possible to extract a convergent
subsequence \( \left\{ b_{k_{l}}\right\} _{l\in \NN } \). Let \( r_{l}=b_{k_{l}} \),
let \( r \) be the rotation of center \( \alpha  \) obtained as
limit \[
r=\lim _{l\toinf }r_{l}=\lim _{l\toinf }b_{k_{l}}\]
and let \( m_{l}=n'_{k_{l}} \); then \begin{eqnarray*}
\lim _{l\toinf }s^{m_{l}}b^{-1}r_{l}s^{-m_{l}} & = & \lim _{l\toinf }s^{n'_{k_{l}}}b^{-1}\left( s^{n'_{k_{l}+i_{k_{l}}}-n'_{k}}bs^{-(n'_{k_{l}+i_{k_{l}}}-n'_{k_{l}})}\right) s^{-n'_{k_{l}}}\\
 & = & \lim _{l\toinf }\left( s^{n'_{k_{l}}}bs^{-n'_{k_{l}}}\right) ^{-1}s^{n'_{k_{l}+i_{k_{l}}}}bs^{-n'_{k_{l}+i_{k_{l}}}}\\
 & = & e.
\end{eqnarray*}

\end{proof}
Using the last lemma along with the uniform continuity properties
that we have when the measure \( \mu  \) is spread out, we are able
to obtain some invariance properties for the limit measure of the
potential kernel; we determine a set of periods for limit measures
\( \nu , \) that is the elements \( \gamma  \) of the group \( \Gamma  \)
such that \[
{\gamma }*\nu =\nu ,\]
 and we show that the limit measure depends only on the speed on of
convergence to the boundary. 

\begin{thm}
\label{prop-Per.Mis.Pot.bT} Suppose that the measure \( \mu  \)
is spread out on \( \Gamma  \). Let \( \left\{ g_{n}\right\} _{n\in \NN } \)
be a sequence in \( \Gamma  \) that converges to \( \alpha \in \bT  \)
and such that the measures \( g_{n}*U \) vaguely converge to \( \nu  \).
Then:

1. There exists a subsequence \( \left\{ g_{n_{k}}\right\} _{k} \)
such that \[
\lim _{k\toinf }{g_{n_{k}}g}*U=\nu \qquad \textrm{for}\, \textrm{all }g\in \Gamma \]

2. Every element of \( \Gamma  \) that fixes \( \alpha  \) is a
period for \( \nu  \).

3. If \( s\in \Gamma  \) is such that \( s\alpha =\alpha  \) and
\( \phi (s)=1 \) then \[
\lim _{n\toinf }g_{n}*U=\lim _{n\toinf }s^{\phi (g_{n})}*U.\]

\end{thm}
\begin{proof}
First suppose that \( g_{k}=s^{n_{k}} \) with \( s\in \Gamma  \)
such that \( s\alpha =\alpha  \) and \( \phi (s)=1 \). 

1. Possibly extracting a subsequence, one knows (Lemma \ref{lem-contuni})
that the Radon measures \[
\nu _{g}=\lim _{k\toinf }{s^{n_{k}}g}*U\]
 are well defined for every \( g\in \Gamma  \) and that this family
depends continuously on \( g \). Therefore the group \[
P=\left\{ b\in \horg :\: \forall g\in \Gamma \: \nu _{bg}=\nu _{g}\right\} \]
is closed. We first prove that \( P \) is normal in \( \Gamma  \).
In fact, if \( \gamma \in \Gamma  \) and if the subsequence \( \left\{ s^{n'_{k}}\gamma s^{-n'_{k}}\right\} _{k} \)
converges to \( \gamma ' \) (see Lemma \ref{lem-conv.alpha}), then
for all \( b\in P \) \begin{eqnarray*}
\nu _{\gamma b\gamma ^{-1}g} & = & \lim _{k\toinf }{s^{n'_{k}}\gamma b\gamma ^{-1}g}*U\\
 & = & \lim _{k\toinf }{s^{n'_{k}}\gamma s^{-n'_{k}}}*{s^{n'_{k}}b\gamma ^{-1}g}*U\\
 & = & {\gamma '}*\lim _{k\toinf }{s^{n'_{k}}b\gamma ^{-1}g}*U\qquad \textrm{using uniform left continuity }\\
 & = & {\gamma '}*\nu _{b\gamma ^{-1}g}\\
 & = & {\gamma '}*\nu _{\gamma ^{-1}g}\\
 & = & {\gamma '}*{\gamma '}^{-1}*\nu _{g}\qquad \textrm{using uniform left continuity }\\
 & = & \nu _{g}
\end{eqnarray*}
so that  \( \gamma b\gamma ^{-1}\in P \).

Next we show that \( P\setminus \horg  \) is compact. In fact the
preceding lemma says that for all \( b\in \horg  \) it is possible
to find a sequence \( \left\{ r_{k}\right\} _{k} \) in \( \horg  \)
that converges to an element \( r\in K_{\alpha } \) and a subsequence
such that \( s^{n'_{k}}b^{-1}r_{k}s^{-n'_{k}} \) converges to the
identity. Then \begin{eqnarray*}
\nu _{b^{-1}rg} & = & \lim _{k\toinf }{s^{n'_{k}}b^{-1}rg}*U\\
 & = & \lim _{k\toinf }{s^{n'_{k}}b^{-1}r_{k}g}*U\qquad \textrm{using uniform right continuity}\\
 & = & \lim _{k\toinf }{s^{n'_{k}}b^{-1}r_{k}s^{-n'_{k}}}*{s^{n'_{k}}g}*U\\
 & = & \lim _{k\toinf }{s^{n'_{k}}g}*U\quad \qquad \textrm{because }s^{n'_{k}}b^{-1}r_{k}s^{-n'_{k}}\rightarrow e\\
 & = & \nu _{g}
\end{eqnarray*}
i.e. \( b^{-1}r\in P \), that is the class \( Pb \) has a representative
\( r \) in \( K_{\alpha } \). Let now \( \pi  \) be the projection
of \( \horg  \) on \( P\setminus \horg  \). Then \( P\setminus \horg =\pi (K_{\alpha }) \)
is compact, because \( \pi  \) is continuous and \( K_{\alpha } \)
is compact.

We will now show that \( \nu _{g}=\nu  \) for every \( g\in \Gamma  \).
Fix a function \( f\in C_{c}(\Gamma ) \) and define the function
\( h \) by \[
h(g)=\nu _{g}(f)\]
Note that \( h \) is bounded because the potential kernel is bounded,
that it is continuous by uniform right continuity and that it is \( \mu  \)-harmonic
on the right because \begin{eqnarray*}
h*\mu (g) & = & \int _{\Gamma }\lim _{k\rightarrow \infty }{s^{n_{k}}g\gamma }*U(f)\mu (d\gamma )\\
 & = & \lim _{k\rightarrow \infty }\int _{\Gamma }{s^{n_{k}}g\gamma }*U(f)\mu (d\gamma )\: \textrm{ by dominated convercence}\\
 & = & \lim _{k\rightarrow \infty }{s^{n_{k}}g}*\sum ^{\infty }_{i=0}\mu ^{(i)}*\mu (f)\\
 & = & \lim _{k\rightarrow \infty }{s^{n_{k}}g}*\sum ^{\infty }_{i=1}\mu ^{(i)}(f)\\
 & = & \lim _{k\rightarrow \infty }{s^{n_{k}}g}*U(f)-f(s^{n_{k}}g)\\
 & = & h(g)\: \textrm{ }\qquad \textrm{because }f\textrm{ has compact support}.
\end{eqnarray*}
 As \( h \) is invariant by left translation of every element of
the group \( P \), it projects onto a function on \( P\setminus \Gamma  \),
\[
\overline{h}(Pg)=h(g),\]
that is harmonic for the measure \( \overline{\mu } \) obtained by
projection of \( \mu  \) on \( P\setminus \Gamma  \). The group
\( P\setminus \Gamma  \) is not Abelian, but taking the left quotient
by the compact subgroup \( P\setminus \horg  \), one obtains \[
\left( P\setminus \horg \right) \setminus \left( P\setminus \Gamma \right) \cong \horg \setminus \Gamma \cong \ZZ .\]
 Thus we are able to use a generalization to the compact extensions
of \( \ZZ  \) of Choquet-Deny theorem, due to Guivarc'h, cf. Théorème
V.2 \cite{Gu73} (an aperiodic measure in this paper is a measure
that generate \( \Gamma  \) as closed group), and conclude that every
continuous bounded harmonic function on \( P\setminus \Gamma  \)
is constant. Hence the function \( \overline{h} \) and therefore
also \( h \) are constant and we can conclude that \[
\nu _{g}(f)=h(g)=h(e)=\nu _{e}(f)=\nu (f)\quad \forall g\in \Gamma .\]

2. We can now show that every \( \gamma \in \Gamma  \) that fixes
\( \alpha  \) is a period of the limit measure \( \nu  \). We first
note that, as \( \gamma  \) fixes \( \alpha  \), the sequence \( s^{-n_{k}}\gamma s^{n_{k}} \)
is relatively compact because \( s^{-n_{k}}\gamma s^{n_{k}}\alpha =\alpha  \)
and \( \phi (s^{-n_{k}}\gamma s^{n_{k}})=\phi (\gamma ) \) (we have
here a situation that is in some sense the opposite to Lemma \ref{lem-conv.alpha},
where we showed that \( s^{n_{k}}\gamma s^{-n_{k}} \) is relatively
compact). Let \( \left\{ s^{-n'_{k}}\gamma s^{n'_{k}}\right\} _{k} \)
be a subsequence that converges to \( \gamma ' \) then by right continuity\[
\gamma *\nu =\lim _{k\toinf }s^{n'_{k}}\left( s^{-n'_{k}}\gamma s^{n'_{k}}\right) *U=\nu _{\gamma '}=\nu .\]

That ends the proof of points 1 and 2 when \( g_{k}=s^{n_{k}} \).

3. Let now \( \left\{ g_{k}\right\} _{k} \) be a generic subsequence
that converges to \( \alpha  \). It can be decomposed in an horocyclic
component and homothetic component by setting \( g_{k}=b_{k}s^{n_{k}} \)
with \( b_{k}\in \horg  \) and \( n_{k}=\phi (g_{k}) \). We observe
that \( n_{k}\ver +\infty  \), so that we are able to apply to \( s^{n_{k}} \)
what we have proved so far. On the other hand the sequence of the
tree ends \( \left\{ b_{k}\alpha =g_{k}\alpha \right\} _{k} \) converges
to \( \alpha  \), thus \( \left\{ b_{k}\right\} _{k} \) is relatively
compact and all its accumulation points are periods of the limit measure
because they fix \( \alpha  \). More precisely for every subsequence
along which \( \left\{ s^{n_{k}}g*U\right\} _{k} \) converges, we
can extract a sub-subsequence along which \( \left\{ b_{k_{l}}\right\} _{l} \)
converges to a rotation \( r \) in \( K_{\alpha } \), then \begin{eqnarray*}
\lim _{l\toinf }{g_{k_{l}}g}*U & = & \lim _{l\toinf }{b_{k_{l}}}*{s^{n_{k_{l}}}g}*U\\
 & = & {r}*\lim _{l\toinf }{s^{n_{k_{l}}}g}*U\\
 & = & \lim _{l\toinf }{s^{n_{k_{l}}}g}*U\quad \textrm{because }r\alpha =\alpha 
\end{eqnarray*}
 Thus \[
\lim _{k\toinf }{g_{k}g}*U=\lim _{k\toinf }{s^{n_{k}}g}*U\]
 and we conclude the proof of point 3. 

The results of point 1 and 2 for a generic sequence \( \left\{ g_{k}\right\} _{k} \)
are a straightforward consequence of the last equality and of the
analogous results for \( g_{k}=s^{n_{k}} \).
\end{proof}
When the group \( \Gamma  \) acts on the tree in a sufficiently homogeneous
and Abelian way, the results of the previous theorem hold even if
we do not assume that \( \mu  \) is spread out. More precisely, we
suppose that \( \Gamma  \) satisfies to the following hypotheses

\begin{lyxlist}{00.00.0000}
\item [(HA)]There exists an end \( \alpha _{0}\in \bT  \) such that the
stabilizer, \( A \), of \( \alpha _{0} \) in \( \Gamma  \) is Abelian.
Besides it exists a measurable set, \( T \), of left coset representatives
 of \( A \) in \( \Gamma  \) \[
\Gamma =TA\]
 and a reference homothety \( c\in A \) with \( \phi (c)=1 \) such
that for every \( t\in T \) we have\[
\lim _{n\toinf }c^{n}tc^{-n}=e\]

\end{lyxlist}
In other words we require that \( \Gamma  \) should be decomposed
as product of an Abelian group of roto-homotheties \( A \) and of
a set of translation \( T \) that acts asymptotically like the identity
far away from the center of the bottom boundary \( \bT  \). This
hypothesis holds if \( \Gamma  \) is contained in \( \affq =\QQ _{p}\rtimes \QQ ^{*}_{p} \)
(in which case we can chose \( T\subseteq \QQ _{p} \) and \( A\subseteq \QQ ^{*}_{p} \)
and the contraction \( c \) may be chosen equal to \( (0,p) \))
and when \( \Gamma  \) is the closure of the Lamplighter group (then
\( T \) is the closure of \( \sZ _{q} \) and \( A=\ZZ  \)).

\begin{prop}
\label{prop-Per.Mis.Pot.bT.Qp}Suppose that \( \Gamma  \) satisfies
to (HA). Let \( \left\{ g_{n}\right\} _{n} \) be a sequence in \( \Gamma  \)
that converges to an end \( \alpha \in \bT  \) and such that \( \left\{ g_{n}*U\right\} _{n} \)
converges to a limit measure \( \nu  \). Then

1. There exists a subsequence \( \left\{ g_{n_{k}}\right\} _{k} \)
such that \[
\lim _{k\toinf }g_{n_{k}}g*U=\nu \qquad \textrm{for}\, \textrm{all }g\in \Gamma \]

2. Every element of \( \Gamma  \) that fixes \( \alpha  \) is a
period for \( \nu  \).

3. If \( s\in \Gamma  \) is such that \( s\alpha =\alpha  \) and
\( \phi (s)=1 \) then \begin{equation}
\label{eq-U(g,.)=U(s,.)Q_p}
\lim _{n\toinf }g_{n}*U=\lim _{n\toinf }s^{\phi (g_{n})}*U.
\end{equation}

\end{prop}
\begin{proof}
The proof follows the same scheme of the previous theorem, being careful
that now the potential kernel is not a priori uniformly right continuous.

Let be \( c \) the contraction that appears in the hypothesis (HA).
First, we suppose that \( \alpha  \) is its center \( \alpha _{0} \),
i.e. \( c\alpha =\alpha  \) and that \( g_{k}=c^{n_{k}} \). For
every function \( f\in C_{c}(\Gamma ) \) and every element \( g=ta \)
of the group \( \Gamma  \), the sequence \( \left\{ c^{n_{k}}g*U(f)\right\} _{k} \)
is bounded and converges to \( a*\nu (f) \). In fact for every convergent
subsequence we have\begin{eqnarray*}
\lim _{k\toinf }c^{n'_{k}}g*U(f) & = & \lim _{k\toinf }c^{n'_{k}}tc^{-n'_{k}}c^{n'_{k}}a*U(f)\\
 & = & \lim _{k\toinf }c^{n'_{k}}a*U(f)\: \quad \textrm{because of the left continuty}\\
 & = & a*\lim _{k\toinf }c^{n'_{k}}*U(f)\quad \textrm{as }A\textrm{ is Abelian}\\
 & = & a*\nu (f)
\end{eqnarray*}
 It is then possible to define the function \( h \) on \( \Gamma  \)\[
h(g)=\lim _{k\toinf }c^{n_{k}}g*U(f)=a*\nu (f),\]
 that is continuous (because \( a*\nu (f) \) is continuous), harmonic
and bounded. As it projects to a continuous bounded function on the
Abelian group \( A \) that is harmonic for the marginal of \( \mu  \)
on \( A \), it is has to be constant by the Choquet-Deny Theorem.
This ends the proof of point 1.  

We have also proved that\[
\nu (f)=a*\nu (f)\quad \textrm{for all }a\in A;\]
 i.e. point 2.

To prove point 3 and to deal with a generic sequence \( \left\{ g_{n}\right\} _{n} \),
we proceed exactly as in the proof of the previous theorem. 

If the sequence \( \left\{ g_{n}\right\} _{n} \) converges to another
end \( \alpha  \), we just need to conjugate by an element of the
group that sends \( \alpha _{0} \) on \( \alpha  \).
\end{proof}

\subsubsection{\textmd{Characterization}}

This section is devoted to the characterization of the of the limit
measures on a neighborhood of the bottom boundary \( \bG  \). This
characterization is given using the decomposition of \( \Gamma  \)
as semi-direct product of \( \ZZ  \) and \( \horg  \). Note that
this decomposition depends on the choice of a reference homothety
\( s \) and that the end \( \alpha \in \bT  \) such that \( s\alpha =\alpha  \)
is then considered as the center of the bottom boundary of the tree.
Also observe that we denote by \( \widehat{\nu } \) the image of
measure \( \nu  \) on a group under the map \( g\mapsto g^{-1} \).

\begin{thm}
\label{teo-lim.U.in.bT} Suppose that \( \esp {\nor{\phi (X_{1})}}<\infty  \).
Then the following holds : 

1. If \( \mu (\phi )=\esp {\phi (X_{1})}>0 \), the only accumulation
point of \( \{g*U\}_{g\in \Gamma } \) when \( g \) converges towards
a point of \( \bG  \) is the null measure. 

If we also suppose that \( \mu  \) is spread out or that \( \Gamma  \)
satisfies to the hypothesis (HA) then

2. If \( \mu (\phi )<0 \) and \( \esp {\nor{X_{1}}}<+\infty  \)
, then for every \( \beta \in \bG  \) and for every \( b\in \Gamma  \)
such that \( b\alpha =\beta  \) \[
\lim _{g\rightarrow \beta }g*U=b*m_{\left\langle s\right\rangle }*\widehat{\overline{m}}\]
where \( m_{\left\langle s\right\rangle } \) is the counting measure
on the subgroup of \( \Gamma  \) generated by \( s \) and \( \overline{m} \)
is the unique \( \hat{\mu } \)-invariant Radon measure on \( \horg  \)
with total mass equal to \( -\frac{1}{\mu (\phi )} \) and invariant
by right action of \( K_{\alpha } \). 

3. If \( \mu (\phi )=0 \), \( \esp {\phi (X_{1})^{2}}<+\infty  \)
and \( \esp {\nor{b(X_{1})}^{2+\varepsilon }}<+\infty  \) then for
every \( \beta \in \bG  \) and for every \( b\in \Gamma  \) such
that \( b\alpha =\beta  \) \[
\lim _{g\rightarrow \beta }g*U=b*m_{\left\langle s\right\rangle }*\widehat{\overline{m}}\]
where \( \overline{m} \) is the unique \( \hat{\mu } \)-invariant
Radon measure on \( \horg  \) defined as extension to \( \horg  \)
of the \( \hat{\mu } \)-invariant measure (\ref{eq-misinv}) on \( \bT  \).
\end{thm}
\begin{proof}
1. If \( \mu (\phi )>0 \), the random walk \( S_{n}=\phi (R_{n}) \)
on \( \ZZ  \) is transient and the only accumulation point of its
potential kernel , \( U_{\phi } \), in a neighborhood of \( +\infty  \)
is zero (cf. proposition 3.4 \cite{Re}). Thus for every bounded non
negative function \( f \) with compact support on \( \Gamma  \)
there exists a bounded non negative function \( F \) with finite
support on \( \ZZ  \) such that \( f(g)\leq F(\phi (g)) \); thus
\[
0\leq \lim _{g\rightarrow \beta }g*U(f)\leq \lim _{g\rightarrow \beta }g*U(F\circ \phi )=\lim _{g\rightarrow \beta }U_{\phi }(\phi (g),F)=0\]
 because \( \phi (g) \) converges to \( +\infty  \) when \( g \)
converges to \( \beta \in \bT  \).

2. and 3. Let \( \left\{ g_{n}\right\} _{n} \) be a sequence that
converges to \( \alpha  \) and such that \( \left\{ g_{n}*U\right\} _{n} \)
converges to a limit measure \( \nu  \). We want to prove that \( \hat{\nu }=m_{\left\langle s\right\rangle }*\widehat{\overline{m}} \)
. 

Using Theorem \ref{prop-Per.Mis.Pot.bT}, we can suppose that \( g_{n}=s^{\phi _{n}} \)
and we know that for all \( g\in \Gamma  \)

\[
\nu =\lim _{n\toinf }{s^{\phi _{n}}g}*U.\]
Let \( \widehat{U} \) be the potential measure associated with the
measure \( \hat{\mu } \), image of \( \mu  \) under group inversion.
Then\[
\widehat{\nu }=\lim _{n\toinf }\widehat{{s^{\phi _{n}}g}*U}=\lim _{n\toinf }\widehat{U}*g^{-1}s^{-\phi _{n}}.\]
 As \( \left\{ \phi _{n}\right\} _{n} \) converges to \( +\infty  \),
for every function \( f \) with compact support in \( \Gamma  \),
the functions \[
x\mapsto \left( f*s^{-\phi _{n}}\right) (x)=f(xs^{-\phi _{n}})\]
 have their support in \( \horg \times \ZZ _{+} \) for sufficiently
large \( n \). As \( \hat{\mu }(\phi )=-\mu (\phi )\geq 0 \), we
can apply to \( \hat{\mu } \) the Corollary \ref{cor-eq.ren.AffT};
thus there exists a probability \( \overline{p} \) on \( \horg  \)
such that for any sufficiently large \( n \) one has

\begin{eqnarray*}
\hat{U}*\overline{p}*s^{-\phi _{n}}(f)=\hat{U}*\overline{p}(f*s^{-\phi _{n}})=\overline{m}*m_{\left\langle s\right\rangle }(f*s^{-\phi _{n}})=\overline{m}*m_{\left\langle s\right\rangle }(f). & 
\end{eqnarray*}
On the other hand, as \( \overline{p} \) has finite mass, by dominated
convergence \[
\lim _{n\toinf }\hat{U}*\overline{p}*s^{-\phi _{n}}(f)=\int _{\Gamma }\lim _{n\toinf }\hat{U}*g*s^{-\phi _{n}}(f)\overline{p}(dg)=\int _{\Gamma }\hat{\nu }(f)\overline{p}(dg)=\hat{\nu }(f).\]
This ends the proof in the case \( \beta =\alpha  \). When \( \beta \neq \alpha  \),
one just need to multiply on the left by an element \( b\in \Gamma  \)
such that \( b\beta =\alpha  \).
\end{proof}

\subsection{Limit near \protect\( \omega \protect \)}

We now study the limit measure in a neighborhood of the mythic ancestor,
\( \omega  \), and show that in this case the limit is always zero. 

\begin{thm}
\label{teo-lim.pot.omega} Suppose that \( \mu  \) is spread out
and that \( \phi (X_{1}) \) is integrable. If \( \mu (\phi )<0 \),
we also suppose that \( \esp {\nor{X_{1}}}<+\infty  \), while if
\( \mu (\phi )=0 \) we suppose that \( \esp {\phi (X_{1})^{2}}<+\infty  \)
and \( \esp {\nor{b(X_{1})}^{2+\varepsilon }}<+\infty  \). Then \[
\lim _{g\ver \omega }g*U=0.\]

\end{thm}
We would like to observe that as, it can be seen in the proof, the
hypothesis that the measure is spread out is not needed when \( \mu (\phi )\neq 0 \)
and \( g \) goes to \( \omega  \) in such a way that \( \phi (g) \)
goes to \( +\infty  \).

\begin{proof}
If \( g \) converges to \( \omega  \) in such a way that \( \phi (g) \)
is bounded from above we can directly apply the Theorem 2.16 in \cite{El82},
which says that on every non-unimodular group, if the probability
law of the random walk is spread out, the potential kernel converges
to zero when \( g \) goes to infinity in such a way that the module
of the Haar measure of the group, that in our case is \( \Delta (g)=q^{\phi (g)} \),
is bounded from above.

Therefore we only need to show that for every sequence \( \left\{ g_{n}\right\} _{n} \)
that converges to \( \omega  \) and such that \( \left\{ \phi (g_{n})\right\} _{n} \)
converges to \( +\infty  \) and for every non-negative continuous
function \( f \) with compact support , \( \left\{ g_{n}*U(f)\right\} _{n} \)
converges to zero. We will distinguish three cases according to sign
of the drift \( \mu (\phi ) \).

Case 1: \( \mu (\phi )>0 \). Exactly as in the proof of Theorem \ref{teo-lim.U.in.bT}.1
in this case one can directly apply the renewal theorem for the induced
random walk on \( \ZZ  \). 

Case 2: \( \mu (\phi )<0 \). First, note that in this case the renewal
theorem on \( \ZZ  \) says that\[
\lim _{h\ver +\infty }U_{\phi }(h,\cdot )=\frac{1}{-\mu (\phi )}m_{\ZZ }\]
where \( m_{\ZZ } \) is the counting measure. On the other side we
have just seen in Theorem \ref{teo-lim.U.in.bT} that if one identifies
\( \ZZ  \) with the subgroup generated by the reference homothety
\( s \) then\[
\lim _{h\ver +\infty }s^{h}*U=m_{\ZZ }*\widehat{\overline{m}}\]
 where \( \widehat{\overline{m}} \) is a measure on \( \horg  \)
whose mass is exactly \( \frac{1}{-\mu (\phi )} \). Then for every
compact set \( H \) in \( \ZZ  \) and every \( \varepsilon >0 \)
there exists a compact open set \( J_{\varepsilon } \) in \( \horg  \)
such that\begin{eqnarray*}
\lim _{h\ver +\infty }s^{h}*U(HJ^{c}_{\varepsilon }) & = & \lim _{h\ver +\infty }s^{h}*U(H\horg -HJ_{\varepsilon })\\
 & = & \lim _{h\ver +\infty }U_{\phi }(h,H)-\lim _{h\ver +\infty }s^{h}*U(HJ_{\varepsilon })\\
 & = & m_{\ZZ }(H)(\frac{1}{-\mu (\phi )}-\widehat{\overline{m}}(J_{\varepsilon }))<\varepsilon ;
\end{eqnarray*}
 i.e. the family of measures \( \left\{ s^{h}*U(H\cdot )\right\} _{h\in \NN } \)
is tight on \( \horg  \). Now fix a compact set \( K \) in \( \Gamma  \)
and observe that  \( g=b(g)s^{\phi (g)} \). Then\[
g*U(K)=U(g^{-1}K)=U(s^{-\phi (g)}b(g)^{-1}K)=s^{\phi (g)}*U(b(g)^{-1}K).\]
 Note that for every sequence \( \left\{ g_{n}\right\} _{n} \) that
converges to \( \omega  \) in such a way that \( \left\{ \phi (g_{n})\right\} _{n} \)
converges to \( +\infty  \), also its projection \( \left\{ b(g_{n})\right\} _{n} \)
on the horocyclic group and its inverse \( \left\{ b(g_{n})\inv \right\} _{n} \)
converge to \( \omega  \). Let \( H \) be a compact set of \( \ZZ  \)
such that \( \phi (K)\subseteq H \), then for every \( \varepsilon >0 \)
and for all \( x\in K \) \[
b(g_{n})\inv x=s^{\phi (x)}s^{-\phi (x)}b(g_{n})\inv s^{\phi (x)}b(x)\in HJ^{c}_{\varepsilon }\qquad \, \, \, \textrm{for any sufficiently large }n\textrm{ }.\]
 Thus, for sufficiently large \( n \), we have \( b(g_{n})\inv K\subseteq HJ^{c}_{\varepsilon } \).
We can conclude that \[
\overline{\lim _{n\ver \infty }}g_{n}*U(K)=\overline{\lim _{n\ver \infty }}s^{\phi (g_{n})}*U(b(g_{n})^{-1}K)\leq \overline{\lim _{n\ver \infty }}s^{\phi (g_{n})}*U(HJ^{c}_{\varepsilon })<\varepsilon \]
i.e. \( g_{n}*U(K) \) converges to zero.

Case 3: \( \mu (\phi )=0 \). Let \( m=\max \{\phi (g):\, g\in \supp \, f\} \)
and for every fixed \( g\in \Gamma  \) let \( t \) be the first
time when \( \phi (g)+S_{n}=\phi (gR_{n}) \) is below \( m \) \[
t=\inf \{k\geq 0:\, \phi (g)+S_{k}\leq m\}.\]
Thus for every \( n<t \) on has \( f(gR_{n})=0 \), and therefore
\[
g*U(f)=\esp {\sum ^{+\infty }_{n=0}f(gR_{n})}=\esp {\sum ^{t-1}_{n=0}f(gR_{n})}+\esp {gR_{t}*U(f)}=\esp {gR_{t}*U(f)}.\]
We have already seen that, when \( \gamma  \) converges to \( \omega  \)
in such a way that \( \phi (\gamma ) \) is bounded from above, \( \gamma *U(f) \)
converges to 0. Hence for every \( \varepsilon >0 \), there is a
compact set \( K_{\varepsilon } \) in \( \Gamma  \) such that \[
\gamma *U(f)\leq \varepsilon \qquad \, \, \textrm{for all }\, \, \gamma \in K^{c}_{\varepsilon }\cap [\phi \leq m].\]
 By definition \( \phi (gR_{t})\leq m \) therefore \begin{eqnarray*}
g*U(f) & = & \esp {(gR_{t}*U(f))\Ind{[\phi (gR_{t})\leq m]}}\\
 & \leq  & \varepsilon +\esp {(gR_{t}*U(f))\Ind{[gR_{t}\in K_{\varepsilon }]}}\\
 & \leq  & \varepsilon +C\esp {\Ind{[gR_{t}\in K_{\varepsilon }]}}
\end{eqnarray*}
where \( C \) is the upper bound of the kernel \( g*U(f) \). Let
now \( \left\{ l^{-}_{k}\right\} _{k} \) the sequence of the ladder
stopping times when \( S_{n}=\phi (R_{n}) \) reaches its minima:
\[
l^{-}_{k}=\min \{n>l^{-}_{k-1}:\, \, S_{n}<S_{l^{-}_{k-1}})\}\, \, \, \textrm{et }\, \, l^{-}_{0}=0.\]
Because \( \phi (R_{t}) \) is strictly smaller then the minimum of
\( \phi (R_{n}) \) for \( n<t \), there exists \( i\in \NN  \)
such that \( t=l^{-}_{i} \) . Let \( U_{l^{-}} \) be the potential
measure of the random walk \( \left\{ R_{l_{k}^{-}}\right\} _{k\in \NN } \)
\begin{eqnarray*}
g*U(f) & \leq  & \varepsilon +C\esp {\Ind{[gR_{t}\in K_{\varepsilon }]}}=\varepsilon +C\esp {\Ind{[gR_{l_{i}^{-}}\in K_{\varepsilon }]}}\\
 & \leq  & \varepsilon +C\esp {\sum _{k=0}^{\infty }\Ind{[gR_{l_{k}^{-}}\in K_{\varepsilon }]}}=\varepsilon +C\, g*U_{l^{-}}(K_{\varepsilon })
\end{eqnarray*}
We remark that the random walk \( R_{l_{k}^{-}} \) satisfies the
hypothesis needed to apply the result of point 2 because: \( \esp {\phi (R_{l_{1}^{-}})}<0 \),
Proposition 4 in \cite{CKW} states that the norm of \( R_{l_{1}^{-}} \)
is integrable and Lemma 2.26 in \cite{El82} guarantees that its law
is spread out (the hypothesis that the group is almost connected that
is assumed there is not necessary to prove this result). Thus \[
\overline{\lim _{g\ver \omega }}g*U(f)\leq \varepsilon +\overline{\lim _{g\ver \omega }}g*U_{l^{-}}(K_{\varepsilon })=\varepsilon ;\]
for every positive \( \varepsilon  \) and we can conclude.
\end{proof}

\section*{Appendix: Proofs of section 2}

\subsection*{Proof of Preposition \ref{prop-munon0}}

\begin{proof}
(1) The distance between \( L_{n}\upsilon  \) and a fixed end \( \alpha  \)
of \( \bT  \) is given by

\[
\Theta (L_{n}\upsilon ,\alpha )=q^{-\phi (L_{n})}\Theta (\upsilon ,L_{n}^{-1}\alpha ).\]
As \( \mu (\phi )<0 \), the random walk \( \phi (L_{n}) \) on \( \ZZ  \)
converges almost surely to \( -\infty  \). On the other hand Theorem
\ref{teofond} states that the right random walk \[
L_{n}^{-1}=X_{1}^{-1}\cdots X_{n}^{-1}=\hat{R}_{n},\]
whose drift is \( \esp {\phi (X_{1}^{-1})}=-\mu (\phi )<0 \), converges
to a random element \( \hat{\xi }_{\infty } \) of   \( \bT  \),
whose law does not charge any point. Thus, for every end \( \upsilon  \)
in \( \bT  \), it exists a sub-set \( \Omega _{\upsilon }\subseteq \Omega  \)
of measure 1 such that on \( \Omega _{\upsilon } \): \[
\lim _{n\rightarrow \infty }\phi (L_{n})=-\infty \quad \textrm{ and }\quad \lim _{n\rightarrow \infty }\hat{R}_{n}=\hat{\xi }_{\infty }\neq \upsilon ;\]
therefore on \( \Omega _{\upsilon } \) \[
\lim _{n\rightarrow \infty }\Theta (L_{n}\upsilon ,\alpha )=+\infty .\]

(2) The probability measure \( m \) is \( \mu  \)-invariant. In
fact, let \( f \) be a bounded continuous function on \( \bT  \)
and let \( X \) be a random variable on \( \affT  \) of law \( \mu  \)
independent from the sequence \( \{X_{n}\}_{n\geq 1} \) then \begin{eqnarray*}
\mu \stackrel{\cdot }{*}m(f) & = & \esp {f(X\xi _{\infty })}=\esp {f(X\lim _{n\rightarrow \infty }X_{1}\cdots X_{n})}\\
 & = & \esp {f(\lim _{n\rightarrow \infty }XX_{1}\cdots X_{n})}\textrm{ as }X\textrm{ acts conituously on }\affT \cup \bT \\
 & = & \esp {f(\xi _{\infty })}=m(f)
\end{eqnarray*}

If \( m' \) is another invariant probability measure and \( \Upsilon _{0} \)
a random variable on \( \bT  \) with law \( m' \), independent from
the increments \( \{X_{n}\}_{n} \), then for every bounded continuous
function \( f \) on \( \bT  \) \begin{eqnarray*}
\eta (f) & = & \esp {f(L_{n}\Upsilon _{0})}=\esp {f(R_{n}\Upsilon _{0})}\\
 & = & \lim _{n\rightarrow \infty }\esp {f(R_{n}\Upsilon _{0})}=\esp {\lim _{n\rightarrow \infty }f(R_{n}\Upsilon _{0})}=\esp {f(\xi _{\infty })}=m(f)
\end{eqnarray*}
by dominated convergence, Theorem \ref{teofond} and since \[
\lim _{g\ver \xi }g\upsilon =\xi \qquad \textrm{for all }\upsilon ,\xi \in \bT .\]
 Thus there is a unique invariant probability measure, and, by the
ergodic theorem, for \( m \)-almost every \( \upsilon  \), the Markov
chain \( L_{n}\upsilon  \) visits infinitely often every set of positive
\( m \)-measure. Furthermore, for open sets one can assure that this
holds for all starting point \( \upsilon  \) (and not only for almost
all) because of the contracting property of the Markov chain:\begin{equation}
\label{eq-cont-forte}
\lim _{n\rightarrow \infty }\Theta (L_{n}\upsilon ,L_{n}\varsigma )=\lim _{n\rightarrow \infty }q^{-\phi (L_{n})}\Theta (\upsilon ,\varsigma )=0.
\end{equation}
 
\end{proof}

\subsection*{Proof of Proposition \ref{prop-Pot.bT=3Dmxl}}

\begin{proof}
For readers convenience, we sketch the proof that formally follows
the same scheme as Proposition 2.1 in \cite{BBE}.

We first show that there exists a probability \( p \) on \( \bT \times \ZZ  \)
such that \begin{equation}
\label{eq-mipotinc}
U_{l}\stackrel{\cdot }{*}p=\frac{1}{\esp {S_{l}}}\left( m_{l}\times \Ind{[0,+\infty [}m_{\ZZ }\right) 
\end{equation}
 where \( \ds U_{l}=\sum _{n=0}^{\infty }\mu _{l}^{(n)} \). Let \( \widetilde{\nu }_{l}=\left( m_{l}\times \Ind{[0,+\infty [}m_{\ZZ }\right)  \).
We observe that for every measurable non-negative function \( f \)
: \begin{eqnarray*}
\mu _{l}\stackrel{\cdot }{*}\widetilde{\nu }_{l}(f) & = & \int _{\bT \times \ZZ }\esp {f(L_{l}\upsilon ,S_{l}+z)\Ind{[z\geq 0]}}m_{l}(d\upsilon )m_{\ZZ }(dz)\\
 & = & \int _{\bT \times \ZZ }\esp {f(L_{l}\upsilon ,z)\Ind{[z-S_{l}\geq 0]}}m_{l}(d\upsilon )m_{\ZZ }(dz)\\
 & \leq  & \int _{\bT \times \ZZ }\esp {f(L_{l}\upsilon ,z)\Ind{[z\geq 0]}}m_{l}(d\upsilon )m_{\ZZ }(dz)\textrm{ }\\
 & = & \widetilde{\nu }_{l}(f)\quad \textrm{because }m_{l}\textrm{ is }\mu _{l}-\textrm{invariant}.\textrm{ }
\end{eqnarray*}
 Thus \[
p':=\widetilde{\nu }_{l}-\mu _{l}\stackrel{\cdot }{*}\widetilde{\nu }_{l}\]
 is a positive measure and one calculates its total mass \( \esp {S_{l}} \).
Furthermore, for every bounded non-negative function \( f=f_{1}\times f_{2} \)
such that \( f_{2} \) has compact support \[
\lim _{n\rightarrow \infty }\mu _{l}^{(n)}\stackrel{\cdot }{*}\widetilde{\nu }_{l}(f)\leq \lim _{n\rightarrow \infty }\left\Vert f_{1}\right\Vert _{\infty }\esp {\int _{S_{l_{n}}}^{+\infty }f_{2}(z)m_{\ZZ }(dz)}=0.\]
Thus \[
U_{l}\stackrel{\cdot }{*}p'(f)=\lim _{n\rightarrow \infty }\sum ^{n}_{k=0}\left( \mu ^{(k)}\stackrel{\cdot }{*}\widetilde{\nu }_{l}-\mu ^{(k)}\stackrel{\cdot }{*}\widetilde{\nu }_{l}\right) (f)=\widetilde{\nu }_{l}\]
 and therefor the probability measure \( \ds p=\frac{p'}{\esp {S_{l}}} \)
verifies (\ref{eq-mipotinc}), on compact sets and thus everywhere. 

Let \( f \) be a non-negative function with support in \( \bT \times \ZZ _{+} \).
To conclude, one has to apply (\ref{eq-mipotinc}) to the non-negative
Borel function \( F(\xi ,z)=\esp {\sum _{0}^{l-1}f(L_{k}\upsilon ,S_{k}+z)} \)
and check that \[
U\stackrel{\cdot }{*}p(f)=U_{l}\stackrel{\cdot }{*}p(F)=\widetilde{\nu }_{l}(F)=m\times m_{\ZZ }(f).\]

\end{proof}
\bibliographystyle{alpha}
\bibliography{ren-tree}

\end{document}